\documentclass[final]{siamltex}

\usepackage{graphicx}

\usepackage{color}
\usepackage{bm}
\usepackage{amsfonts}
\usepackage[tight,footnotesize]{subfigure}

\def\norm#1{\|#1\|}
\newcommand{\tr}{^{\sf T}}
\newcommand{\m}[1]{{\bf{#1}}}
\newcommand{\g}[1]{\bm #1}
\newcommand{\C}[1]{{\cal {#1}}}
\newcommand{\prox}{{\rm prox}}
\newtheorem{remark}{Remark}[section]

\title{Inexact alternating direction multiplier methods for
separable convex optimization
\thanks{April 8, 2016.
The authors gratefully acknowledge support by the National
Science Foundation under grants 1522629 and 1522654, and
by the Office of Naval Research under grant N00014-15-1-2048.
}}
\author{
    William W. Hager\thanks{{\tt hager@math.ufl.edu},
        http://people.clas.ufl.edu/hager/,
        PO Box 118105,
        Department of Mathematics,
        University of Florida, Gainesville, FL 32611-8105.
        Phone (352) 294-2308. Fax (352) 392-8357.}
\and
    Hongchao Zhang\thanks{{\tt hozhang@math.lsu.edu},
        http://www.math.lsu.edu/$\sim$hozhang,
        Department of Mathematics,
        Louisiana State University, Baton Rouge, LA 70803-4918.
        Phone (225) 578-1982. Fax (225) 578-4276.}
}
\begin{document}
\maketitle
\begin{abstract}
Inexact alternating direction multiplier methods (ADMMs) are
developed for solving general separable convex optimization problems
with a linear constraint and with an objective that is the sum of smooth
and nonsmooth terms.
The approach involves linearized subproblems, a back substitution step,
and either gradient or accelerated gradient techniques.
Global convergence is established.
The methods are particularly useful when the ADMM subproblems
do not have closed form solution or when the solution of the subproblems
is expensive.
Numerical experiments based on image reconstruction problems
show the effectiveness of the proposed methods.

\end{abstract}


\pagestyle{myheadings}
\thispagestyle{plain}
\markboth{W. W. HAGER AND H. ZHANG}
{INEXACT ADMM FOR SEPARABLE CONVEX OPTIMIZATION}
\section{Introduction}

We consider a convex separable linearly constrained optimization problem
\begin{equation}\label{Prob}
\min \; \Phi (\m{x}) \; \mbox{ subject to } \m{Ax} = \m{b}
\end{equation}
where $\Phi : \mathbb{R}^n \rightarrow \mathbb{R}\cup\{\infty\}$ and
$\m{A}$ is $N$ by $n$.
By a separable convex problem,
we mean that the objective function is a sum of $m$ independent parts,
and the matrix is partitioned compatibly as in
\begin{equation}\label{ProbM}
\Phi(\m{x}) = \sum_{i=1}^m f_i(\m{x}_i) + h_i(\m{x}_i)
\quad \mbox{and} \quad \m{Ax} = \sum_{i =1}^m \m{A}_i \m{x}_i.
\end{equation}
Here $f_i$ is convex and Lipschitz continuously differentiable,
$h_i$ is a proper closed convex function (possibly nonsmooth),
$\m{A}_i$ is $N$ by $n_i$ with $\sum_{i=1}^m n_i = n$,
and the columns of $\m{A}_i$ are linearly independent for $i \ge 2$.
Constraints of the form $\m{x}_i \in \C{X}_i$, where $\C{X}_i$ is a closed
convex set, can be incorporated in the optimization problem by
setting $h_i(\m{x}_i) = \infty$ when $\m{x}_k \not\in \C{X}_i$.
The problem (\ref{Prob})--(\ref{ProbM}) has attracted extensive research 
due to its importance in areas such as
image processing, statistical learning and compressed sensing.
See the recent survey \cite{Boyd10} and its references.

Let $\C{L}$ be the Lagrangian given by
\[
\C{L}(\m{x}, \g{\lambda}) = \Phi(\m{x}) +
\langle \g{\lambda}, \m{Ax} - \m{b} \rangle ,
\]
where $\g{\lambda}$ is the Lagrange multiplier for the linear constraint
and $\langle \cdot , \cdot \rangle$ denotes the Euclidean inner product.
It is assumed that there exists a solution $\m{x}^*$ to
(\ref{Prob})--(\ref{ProbM}) and an associated Lagrange multiplier
$\g{\lambda}^* \in \mathbb{R}^N$ such that $\C{L}(\cdot, \g{\lambda}^*)$
attains a minimum at $\m{x}^*$, or equivalently, the following
first-order optimality conditions hold:
$\m{Ax}^* = \m{b}$ and
for $i =$ $1, 2, \ldots, m$ and for all $\m{u} \in \mathbb{R}^{n_i}$, we have
\begin{equation}\label{Wstar}
\langle \nabla f_i (\m{x}_i^*) + \m{A}_i \tr \g{\lambda}^*,
\m{u} - \m{x}_i^* \rangle + h_i(\m{u}) \ge h_i (\m{x}_i^*),
\end{equation}
where $\nabla$ denotes the gradient.

A popular strategy for solving (\ref{Prob})--(\ref{ProbM})
is the alternating direction multiplier method (ADMM)
\cite{GM76, gl84} given by
\begin{equation}\label{adm}
\quad \quad \left\{
\begin{array}{lcl}
\m{x}_i^{k+1} &=&
\arg \displaystyle{\min_{\m{x}_i \in \mathbb{R}^{n_i}}}
\; L (\m{x}_1^{k+1}, \ldots,\m{x}_{i-1}^{k+1},\m{x}_i,
\m{x}_{i+1}^k, \ldots, \m{x}_m^k , \g{\lambda}^k ), \; i=1, \ldots, m, \\
\g{\lambda}^{k+1} &=&  \g{\lambda}^k + \rho
(\m{Ax}^{k+1} - \m{b}), 
\end{array}
\right.
\end{equation}
where $L$, the augmented Lagrangian, is defined by
\begin{equation}\label{AL}
L (\m{x}, \g{\lambda}) = \C{L}(\m{x}, \g{\lambda}) +
\frac{\rho}{2} \| \m{Ax} - \m{b} \|^2.
\end{equation}
Here $\rho >0$ is the penalty parameter.
Early ADMMs only consider problem (\ref{Prob})--(\ref{ProbM})
with $m=2$ corresponding to a $2$-block structure.
In this case, the global convergence and complexity
can be found in \cite{EB92, HeYuan12}.
When $m \ge 3$ the ADMM strategy (\ref{adm}),
a natural extension of the $2$-block ADMM,
is not necessarily convergent \cite{chyy2016},
although its practical efficiency has been observed in 
many recent applications \cite{TaoYuan2011,wgy2010}.

Recently, much research has focused on modifications to ADMM to ensure
convergence when $m \ge 3$.
References include \cite{CaiHanYuan14,ChenLiLiuYe15,ChenShenYou13,
DavisYin15, GoldfarbMa2012, HanYuan12, HTXY2013, HeTaoXuYuan12,
LiSunToh14, LinMaZhang14, LinMaZhang15}.
One approach \cite{CaiHanYuan14, ChenShenYou13, DavisYin15, HanYuan12}
assumes $m-2$ of the functions in the objective are strongly convex
and the penalty parameter is sufficiently small.
Linear convergence results under additional conditions are obtained
in \cite{LinMaZhang15}.
Analysis of a randomly permuted ADMM which allows for nonseparated variables is
given in \cite{ChenLiLiuYe15}.
Another approach, first developed in \cite{HTXY2013, HeTaoXuYuan12},
involves a back substitution step to complement
the ADMM forward substitution step.
The algorithms developed in our paper utilize this back
substitution step.

The dominant computation in an iteration of ADMM is 
the solution of the subproblems in (\ref{adm}).
The efficiency depends on our
ability to solve these subproblems inexactly while 
still maintaining the global convergence,
especially when no closed formula exists for the
subproblems \cite{YangZhang2011}.
One line of research is to solve the subproblems to an accuracy
based on an absolute summable error criterion
\cite{ChenSunToh2015, EB92, GolTre1979}.
In \cite{LiLiaoYuan2013}, the authors combine an adaptive error criterion
with the absolute summable error criterion for 2-block ADMM with
logarithmic-quadratic proximal regularization and
further correction steps to modify the solutions generated
from the ADMM subproblems.
In \cite{EcksteinYao16}, the authors develop a 2-block ADMM with a relative
error stopping condition for the subproblems, motivated by \cite{Eckstein13},
based on the total subgradient error.
Another line of research is to add proximal terms to make the subproblems
strongly convex \cite{ChenTeboulle94,HLHY2002} and relatively easy to solve.
However, this approach often requires accurate solution of
the proximal subproblems.
When $m=1$, ADMM reduces to the standard augmented Lagrangian method (ALM), 
for which practical relative error criteria for solving 
the subproblems have been developed and encouraging numerical results
have been obtained \cite{Eckstein13, SolSva2000}.
In this paper, motivated by our recent work on variable stepsize
Bregman operator splitting methods (BOSVS),
by recent complexity results for gradient 
and accelerated methods for convex optimization, and by the adaptive
relative error strategy used in ALM, we develop new
{\it inexact} approaches for solving the ADMM subproblems.
To the best of our knowledge,
these are the first ADMMs for solving the general separable
convex optimization problem (\ref{Prob})--(\ref{ProbM}) based on an
adaptive accuracy condition that does not employ
an absolute summable error criterion and that guarantees
global convergence, even when $m \ge 3$. 

To guarantee global convergence, a block Gaussian backward substitution
strategy is used to make corrections to the approximate subproblem solutions.
In the special case $m = 2$, 
the method will reduces to a 2-block ADMM without back substitution.
This idea of using  block Gaussian back substitution was first proposed
in \cite{HTXY2013,HeTaoXuYuan12}. 
The method in this earlier work
requires the exact solution of the subproblems to obtain global convergence, 
while our new approach allows an inexact solution.
More recently, a linearly convergent ADMM was developed in    
\cite{HongLuo2013}.
This algorithm linearizes the subproblems to achieve an inexact solution,
and requires that the functions $f_i$ and $h_i$ in the objective
function satisfy certain ``local error bound'' conditions.
In addition, to  ensure linear convergence,
the stepsize $\alpha_k$ in (\ref{adm}) must be sufficiently small,
which could significantly deteriorate the practical performance. 

In this paper, we focus on problems where the minimization of the
dual function over one or more of the primal variable $\m{x}_i$ is nontrial,
and the accuracy of an inexact minimizer needs to be taken into account.
On the other hand, when these minimizations are simple enough,
it is practical to minimize the Lagrangian over $\m{x}$ and
compute the dual function.
This leads to a possibly nonsmooth dual function which can be approached
through smoothing techniques as in \cite{LiChenDongWu16}, or through
active set techniques as in \cite{hz15}.

Our paper is organized as follows.
In the Section \ref{BOSVS},
we  first generalize the  BOSVS algorithm
\cite{chy13, hyz14} to handle multiple blocks.
The original BOSVS algorithm was tailored to the two block case,
but used an adaptive stepsize when solving the subproblem,
and consequently, it achieved much better overall efficiency
when compared to the Bregman operator splitting (BOS) type algorithms
based on a fixed smaller stepsize.
In Sections~\ref{mBOSVS} and \ref{aBOSVS}, more adaptive stopping criteria
for the subproblems are proposed. 
The adaptive criteria for bounding the accuracy in the ADMM subproblems are
based on both the current and accumulated iteration change in the subproblem.
These novel stopping criteria are motivated by the complexity analysis
of gradient methods for convex optimization, and by the relative accuracy
strategy often used in an inexact augmented Lagrangian method
for nonlinear programming.
Although our analysis is carried out 
with vector variables, these results could be extended to
matrix variables which could have more potential applications.

\subsection{Notation}
The set of solution/multiplier pairs
for (\ref{Prob}) is denoted $\C{W}^*$, while
$(\m{x}^*, \g{\lambda}^*) \in \C{W}^*$ is
a generic solution/multiplier pair.
For $\m{x}$ and $\m{y} \in \mathbb{R}^n$,
$\langle \m{x}, \m{y} \rangle = \m{x} \tr \m{y}$ is the standard inner product,
where the superscript $\tr$ denotes transpose.
The Euclidean norm, denoted $\norm{\cdot}$, is defined by
$\norm{\m{x}} =\sqrt{\langle \m{x}, \m{x} \rangle}$ and 
$\norm{\m{x}}_\m{G} =\sqrt{\m{x} \tr \m{G} \m{x}}$
for a positive definite matrix $\m{G}$.
$\mathbb{R}^+$ denotes the set of nonnegative real numbers,
while $\mathbb{R}^{++}$ denotes the set of positive real numbers.
For a differentiable function $f: \mathbb{R}^n \to \mathbb{R}$,
$\nabla f (\m{x})$ is the gradient of $f$ at $\m{x}$, a column vector.
More generally, $\partial f (\m{x})$ denotes the subdifferential at $\m{x}$.
If $\m{x}$ is a vector, then $\m{x}_+$ denotes the subvector obtained
by dropping the first block of variables from $\m{x}$.
Thus if $\m{x} \in \mathbb{R}^n$ with $\m{x}_i \in \mathbb{R}^{n_i}$
for $i \in [1, m]$, then $\m{x}_+ =$
$(\m{x}_2, \m{x}_3, \ldots, \m{x}_m)\tr$.
\section{Algorithm Structure}
\label{algorithm}
Three related inexact ADMMs are developed called generalized,
multistep, and accelerated BOSVS.
They differ in the details of the formula for
the new iterate $\m{x}^{k+1}$, but the overall structure of the
algorithms is the same.
Both multistep and accelerated BOSVS typically represent
a more exact ADMM iteration when compared to generalized BOSVS,
while the accelerated BOSVS subiterations often converge more
rapidly than those of multistep BOSVS.
The common elements of these
three algorithms appear in Algorithm~\ref{ADMMcommon}.
\renewcommand\figurename{Alg.}
\begin{figure}[h]
{\tt
\begin{tabular}{l}
\hline
{\bf Parameters:}
$\rho, \, \sigma, \, \delta_{\min}, \, \theta_1, \, \theta_2,\,
\theta_3 \in \mathbb{R}^{++}$,
$\delta_{\min} <  \delta_{\max}$,
$\alpha \in (0, 1)$, $\sigma < 1 < \tau \le \eta$. \\[.05in]
{\bf Starting guess:} $\m{x}^1$ and $\g{\lambda}^1$. \\[.05in]
{\bf Initialize:} $\m{y}^1 = \m{x}^1$, $k = 1$,
$\delta_{\min, i} = \delta_{\min}$ and
$\Gamma_i^0 = 0$, $1 \le i \le m$, $e^0 = \infty$ \\[.05in]
\begin{tabular}{ll}
{\bf Step 1:} &  For $i=1, \dots, m$  \\[.05in]
& $\quad$ Generate $\m{x}_i^{k+1}$ and $\m{z}_i^k$,  estimate
$r_i^k \approx \|\m{x}_i^{k+1} - \m{x}_i^k\|^2$. \\
& End \\[.05in]
{\bf Step 2:}  & If
$e^k = \theta_1 \|\m{z}_+^{k} - \m{y}_+^k\| + \theta_2 
\| \m{Az}^{k} - \m{b}\| + \theta_3 \sqrt{\sum_{i =1}^m {r}_i^k}$ \\
& is sufficiently small, terminate. \\[.1in]
{\bf Step 3:} & Set
$\m{y}_1^{k+1} = \m{z}_1^{k}$,
$\m{y}_+^{k+1} =$
$\m{y}_+^k + \alpha \m{M}^{-\sf T} \m{H} (\m{z}_+^{k} - \m{y}_+^k) $,
\\[.05in]
& $\g{\lambda}^{k+1} = \g{\lambda}^k + \alpha \rho ( \m{Az}^{k} - \m{b} )$,
$k:=k+1$, and go to Step 1.  \\
\end{tabular}\\
\hline
\end{tabular}
}
\caption{\rm Our ADMM structure.}
\label{ADMMcommon} 
\end{figure}
\renewcommand\figurename{Fig.}

The algorithms generate three sequences
$\m{x}^k$, $\m{y}^k$, and $\m{z}^k$.
In Step~1 of Algorithm~\ref{ADMMcommon} there may be more than one
ADMM subiteration, as determined by an adaptive stopping criterion.
The iterate $\m{x}^{k+1}$ is the final iterate generated in the
ADMM (forward substitution) subproblems,
$\m{y}^k$ is generated by the back substitution process in Step~3, and
$\m{z}^k$ is an average of the iterates in the ADMM subproblems of Step~1.
In generalized BOSVS, $\m{z}^k = \m{x}^{k+1}$ since there is only one
ADMM subiteration, while multistep and accelerated BOSVS typically
perform more than one subiteration and $\m{z}^k$ is obtained
by a nontrivial averaging process.
The matrix $\m{M}$ in Step~3 is the $m-1$ by $m-1$
block lower triangular matrix defined by
\begin{equation}\label{m-def}
\m{M}_{ij} = \left\{ \begin{array}{cl}
\m{A}_{i+1}\tr \m{A}_{j+1} & \mbox{if } 1 \le j \le i < m, \\
\m{0}              & \mbox{if } 1 \le i < j < m.
\end{array} \right.
\end{equation}
The matrix $\m{H}$ is the $m-1$ by $m-1$ block diagonal matrix whose
diagonal blocks match those of $\m{M}$.
The matrices $\m{M}$ and $\m{H}$ are invertible since the columns of
$\m{A}_i$ are linearly independent for $i \ge 2$, which implies that
$\m{A}_i\tr\m{A}_i$ is invertible for $i \ge 2$.

\section{Generalized BOSVS}
\label{BOSVS}
Our first algorithm is a generalization of the BOSVS algorithm developed in
\cite{chy13, hyz14} for a two-block optimization problem.
Let $\Phi_i^k : \mathbb{R}^{n_i} \times \mathbb{R}^{n_i} \times \mathbb{R}$
$\rightarrow \mathbb{R}$ be defined by
\[
\Phi_i^k (\m{u}, \m{v}, \delta) =
f_i(\m{v}) +
\langle \nabla f_i(\m{v}), \m{u} - \m{v} \rangle +  
\frac{\delta}{2} \| \m{u} - \m{v} \|^2 + h_i(\m{u}) +
\frac{\rho}{2} \|\m{A}_i \m{u} - \m{b}_i^k + \g{\lambda}^k/\rho \|^2,
\]
where
\begin{equation}\label{bik}
\m{b}_i^k = \m{b} - \sum_{j < i} \m{A}_j \m{z}_j^{k} -
\sum_{j > i} \m{A}_j \m{y}_j^k .
\end{equation}
The function $\Phi_i^k$ corresponds to the part of the augmented Lagrangian
associated with the $i$-th component of $\m{x}$, but with
the smooth term $f_i$ linearized around $\m{v}$ and with a proximal term
added to the objective.
Algorithm~\ref{1} which follows is the Step~1 inner loop of
Algorithm~\ref{ADMMcommon} for generalized BOSVS.
Throughout the paper, the generalized BOSVS algorithm refers
to Algorithm~\ref{ADMMcommon} with the Step~1 inner loop
given by Algorithm~\ref{1}.

Although $\m{y}^k$ does not appear explicitly in Algorithm~\ref{1},
it is hidden inside the $\m{b}_i^k$ term of $\Phi_i^k$.
The iterate $\m{x}_i^{k+1}$ is obtained by minimizing the
$\Phi_i^k$ function and checking the line search condition of Step~1b.
In Step~1a, mid denotes median and
the initial stepsize $\delta_{i,0}^k$ of Step~1a is a safeguarded version
of the Barzilai-Borwein formula \cite{bb88}.
\renewcommand\figurename{Alg.}
\begin{figure}[h]
{\tt
\begin{tabular}{l}
\hline
{\bf Inner loop of Step 1:}\\[.05in]
\begin{tabular}{ll}
1a. & For $k > 1$, $\delta_{i,0}^k = \mbox{mid }
\left\{ \delta_{\min,i}, \; s^{BB},\; \delta_{\max} \right\}$
and\\[.05in]
& $\quad s^{BB} =$
$\langle \nabla f_i (\m{x}_i^k) - \nabla f_i (\m{x}_i^{k-1}), \;
\m{x}_i^k - \m{x}_i^{k-1}\rangle /\|\m{x}_i^k - \m{x}_i^{k-1}\|^2.$\\[.05in]
& For $k = 1$, $\delta_{i,0}^k$ can be any scalar in
$[\delta_{\min}, \delta_{\max}]$.\\[.05in]
1b. & Set $\delta_i^k = \eta^j \delta_{i,0}^k$, where
$j \ge 0$ is the smallest integer such that \\[.05in]
& $\quad$ $f_i(\m{x}_i^k) +$
$\langle \nabla f_i(\m{x}_i^k), \m{x}_i^{k+1} -\m{x}_i^k \rangle
+ \frac{(1-\sigma)\delta_i^k}{2} \|\m{x}_i^{k+1} - \m{x}_i^k\|^2
\ge f_i(\m{x}_i^{k+1})$,\\[.05in]
& where $\m{x}_i^{k+1} = \m{z}_i^k = \arg \min \{
\Phi_i^k (\m{u}, \m{x}_i^k, \delta_i^k) : \m{u} \in \mathbb{R}^{n_i} \}$.
\\[.05in]
1c. & Set ${r}_i^k = (1/\delta_i^k) \| \m{x}_i^{k+1} -\m{x}_i^k \|^2$.\\[.05in]
1d. & If $k > 1$ and $ \delta_i^k > \max \{\delta_i^{k-1}, \delta_{\min,i} \}$,
then $\delta_{\min,i} := \tau \delta_{\min,i} $. \\
\end{tabular}\\
\hline
\end{tabular}
}
\caption{Inner loop in Step~$1$ of Algorithm~$\ref{ADMMcommon}$
for the generalized BOSVS scheme.}
\label{1} 
\end{figure}
\renewcommand\figurename{Fig.}

Let $\zeta_i$ denote the Lipschitz constant for $\nabla f_i$.
By a Taylor expansion of $f_i$ around $\m{x}_i^k$,
we see that the line search condition of Step~1b is satisfied whenever
$(1-\sigma)\delta_i^k \ge \zeta_i$, or equivalently, when
\begin{equation}\label{xx}
\delta_i^k \ge \zeta_i/(1-\sigma).
\end{equation}
Since $\eta > 1$, $\delta_i^k$ increases as $j$ increases, and consequently,
(\ref{xx}) holds for $j$ sufficiently large.
Hence, if $j > 0$ at the termination of the line search, we have
$\delta_i^k \le$ $\eta \zeta_i/(1-\sigma)$.
If the line search terminates for $j = 0$, then
$\delta_i^k \le$ $\delta_{\max}$.
In summary, we have
\begin{equation}\label{delta_bound}
\delta_{\min} \le \delta_i^k \le
\max \{\eta \zeta_i/(1-\sigma), \delta_{\max} \} \quad
\mbox{for all } k.
\end{equation}

Since $s^{BB} \le \zeta_i$ in Step~1a, it follows that
$\delta_{i,0}^k =$ $\delta_{\min, i}$ whenever
$\delta_{\min, i} \ge \zeta_i$.
In Step~1d, $\delta_{\min, i}$
is increased by the factor $\tau$
whenever $\delta_i^k > \delta_i^{k-1}$.
Hence, after a finite number of iterations where
$\delta_i^k > \delta_i^{k-1}$, we have
$\delta_{\min, i} \ge \zeta_i/(1-\sigma)$, which implies that
the line search terminates at $j = 0$ with
$\delta_i^k =$ $\delta_{i,0}^k =$ $\delta_{\min, i}$.
We conclude that
\begin{equation}\label{delta_monotone}
\delta_i^{k} \le \delta_i^{k-1} \quad \mbox{for $k$ sufficiently large},
\end{equation}
where $\delta_i^k$ denotes the final accepted value in Step~1b.
Note that the inequality in (\ref{delta_monotone}) cannot be replaced
by equality since the number of iterations where
$\delta_i^k > \delta_i^{k-1}$ may not be enough to yield
$\delta_{\min, i} \ge \zeta_i/(1-\sigma)$.

In the BOS algorithm, the line search is essentially eliminated by taking
$\delta_i^k$ larger than the Lipschitz constant $\zeta_i$.
Taking $\delta_i^k$ large, however, causes $\|\m{x}_i^{k+1} - \m{x}_i^k\|$
to be small due to the proximal term in the objective
$\Phi_i^k(\;\cdot\;, \m{x}_i^k, \delta_i^k)$ associated with $\m{x}_i^{k+1}$.
These small steps lead to slower convergence than what is achieved
with BOSVS where $\delta_i^k$ is adjusted by the line search criterion
in order to achieve a small, but acceptable, choice for $\delta_i^k$.

In our analysis of generalized BOSVS, we first observe that when
$e^k =0$, we have reached a solution of (\ref{Prob})--(\ref{ProbM}).
\smallskip

\begin{lemma}\label{L-stop-cond}
If $e^{k}=0$ in the generalized BOSVS algorithm,
then $\m{x}^{k+1} = \m{x}^k = \m{y}^k$ solves $(\ref{Prob})$--$(\ref{ProbM})$
and $(\m{x}^k, \g{\lambda}^k) \in \C{W}^*$.
\end{lemma}
\smallskip

\begin{proof}
Let $\m{x}^*$ denote $\m{x}^k$.
If $e^k=0$, then $r_i = 0$ for each $i$, and by Step~1c of generalized BOSVS,
$\m{x}^{k+1} = \m{x}^k = \m{x}^*$.
For generalized BOSVS, we set $\m{z}^k = \m{x}^{k+1}$ in Step~1b.
Since $\m{x}^{k+1} = \m{x}^k$, it follows that $\m{z}^k =$ $\m{x}^*$.
The identity $\m{z}^k = \m{x}^{k+1}$ also implies that
$\m{z}^{k-1} =$ $\m{x}^k =$ $\m{x}^{*}$.
In Step~3 of Algorithm~\ref{ADMMcommon}, 
$\m{y}_1^k = \m{z}_1^{k-1} = \m{x}_1^*$.
Since $e^k = 0$, Step~2 of Algorithm~\ref{ADMMcommon} implies that
$\m{y}_+^k =$ $\m{z}_+^k =$ $\m{x}_+^*$.
Hence, $\m{y}^k = \m{x}^* = \m{z}^k$.
Since $e^k = 0$, it also follows from Step~2 that $\m{Ax}^* = \m{b}$.
Consequently, we have
\begin{equation}\label{BIK}
\m{b}_i^k = \m{b} - \sum_{j < i} \m{A}_j \m{z}_j^{k} -
\sum_{j > i} \m{A}_j \m{y}_j^k =
\m{b} - \sum_{j < i} \m{A}_j \m{x}_j^{*} -
\sum_{j > i} \m{A}_j \m{x}_j^* = \m{A}_i \m{x}_i^* .
\end{equation}
By Step~1b, $\m{x}_i^{k+1}$ is the minimizer of
$\Phi_i^k(\; \cdot \; , \m{x}_i^k, \delta_i^k)$.
Since $\m{x}_i^{k+1} = \m{x}_i^k = \m{x}_i^*$, it follows that
$\m{x}_i^*$ is the minimizer of
$\Phi_i^k(\; \cdot \; , \m{x}_i^*, \delta_i^k)$.
After taking into account (\ref{BIK}),
the first-order optimality condition associated with the minimizer
$\m{x}_i^*$ of $\Phi_i^k(\; \cdot \; , \m{x}_i^*, \delta_i^k)$
is exactly the same as (\ref{Wstar}), but with $\g{\lambda}^*$ replaced
$\g{\lambda}^k$.
Hence, $(\m{x}^*, \g{\lambda}^k) \in \C{W}^*$.
\end{proof} 

Two lemmas are needed for the convergence of the generalized BOSVS algorithm.
\smallskip

\begin{lemma}\label{lem-prop1}
Suppose that $\m{u}$ and $\m{v} \in \mathbb{R}^{n_i}$ satisfy
\begin{equation}\label{yy}
f_i(\m{v}) + \langle \nabla f_i(\m{v}), \m{u} -\m{v} \rangle
+ \frac{(1-\sigma)\delta}{2} \|\m{u}- \m{v}\|^2 \ge f_i(\m{u})
\end{equation}
for some $\sigma \in [0,1)$ and $\delta >0$,
where $\m{u}$ minimizes $\Phi_i^k(\;\cdot\;,\m{v}, \delta)$.
Then, for any $\m{w} \in \mathbb{R}^{n_i}$ we have
\begin{equation}\label{prop1}
\quad \quad L_i^k (\m{w})  - L_i^k (\m{u})
\ge \frac{\delta}{2} (\|\m{w} - \m{u}\|^2 - \|\m{w} - \m{v}\|^2 )
+ \frac{\rho}{2} \| \m{A}_i (\m{w} - \m{u})\|^2
+ \frac{\sigma\delta}{2}  \|\m{u}- \m{v}\|^2,
\end{equation}
where $L_i^k$ is given by
\begin{equation}\label{lik}
L_i^k(\m{u}) = f_i (\m{u}) +
h_i (\m{u}) +
\frac{\rho}{2} \|\m{A}_i \m{u} - \m{b}_i^k + \g{\lambda}^k/\rho\|^2.
\end{equation}
\end{lemma}
\smallskip

\begin{proof}
Adding $ h_i (\m{u}) +
\frac{\rho}{2} \|\m{A}_i \m{u} - \m{b}_i^k + \g{\lambda}^k/\rho\|^2$ to
each side of the inequality (\ref{yy}) and rearranging, we obtain
\[
\Phi_i^k(\m{u}, \m{v}, \delta) -
\frac{\sigma \delta}{2} \|\m{u}- \m{v}\|^2 \ge 
L_i^k(\m{u}) .
\]
Adding $L_i^k (\m{w})$ to each side of this inequality gives
\begin{equation}\label{zz}
L_i^k (\m{w})  - L_i^k (\m{u})
\ge L_i^k (\m{w}) - \Phi_i^k(\m{u}, \m{v}, \delta) 
+ \frac{\sigma \delta}{2} \|\m{u}- \m{v}\|^2 .
\end{equation}
Utilizing the convexity inequality
$f_i (\m{w}) - f_i(\m{v}) \ge \langle \nabla f_i(\m{v}), \m{w}-\m{v} \rangle$,
we have
\begin{eqnarray*}
L_i^k (\m{w}) - \Phi_i^k(\m{u}, \m{v}, \delta) &\ge&
\frac{\rho}{2}
\left( \|\m{A}_i \m{w} - \m{b}_i^k - \g{\lambda}^k/\rho\|^2 -
\|\m{A}_i \m{u} - \m{b}_i^k - \g{\lambda}^k/\rho\|^2
\right) \\
&& \quad \quad + \langle \nabla f_i(\m{v}), \m{w}-\m{u} \rangle
- \frac{\delta}{2}\|\m{u} - \m{v}\|^2 + h_i(\m{w}) - h_i(\m{u}) .
\end{eqnarray*}
Expand the smooth terms involving $\m{w}$ on the right side
in a Taylor series around $\m{u}$ to obtain
\[
L_i^k (\m{w}) - \Phi_i^k(\m{u}, \m{v}, \delta) \ge
\left\langle \m{g}_i^k, \; \m{w} - \m{u} \right\rangle
+h_i(\m{w}) - h_i (\m{u})
+ \frac{\rho}{2} \|\m{A}_i(\m{w} - \m{u})\|^2
- \frac{\delta}{2} \|\m{u} - \m{v}\|^2 ,
\]
where $\m{g}_i^k =$
$\nabla f_i (\m{v})  + \rho \m{A}_i  \tr (\m{A}_i \m{u} -
\m{b}_i^k+ \g{\lambda}^k/\rho)$.
Since $\Phi_i^k$ is the sum of smooth and nonsmooth terms,
the first-order optimality condition for the minimizer $\m{u}$ of
$\Phi_i^k$ can be expressed
\[
\left\langle \m{g}_i^k + \delta (\m{u} - \m{v}), \; \m{w} - \m{u} \right\rangle
+ h_i(\m{w}) \ge h_i (\m{u})
\]
for all $\m{w} \in \mathbb{R}^{n_i}$, which implies that
\[
\langle \m{g}_i^k, \; \m{w} - \m{u} \rangle
+ h_i(\m{w}) - h_i (\m{u}) \ge - \delta \langle \m{u} - \m{v}, \m{w} - \m{u}
\rangle
\]
for all $\m{w} \in \mathbb{R}^{n_i}$.
Utilizing this inequality, we have
\[
L_i^k (\m{w}) - \Phi_i^k(\m{u}, \m{v}, \delta) \ge
-\delta \left\langle \m{u} - \m{v}, \m{w} - \m{u} \right\rangle
+ \frac{\rho}{2} \|\m{A}_i(\m{w} - \m{u})\|^2
- \frac{\delta}{2} \|\m{u} - \m{v}\|^2 .
\]
Insert this in (\ref{zz}).
Since
\[
2\langle \m{u} - \m{v}, \m{w} - \m{u} \rangle
+\|\m{u}- \m{v}\|^2 =
\|\m{w} - \m{v}\|^2 - \|\m{w} - \m{u}\|^2,
\]
the proof is complete.
\end{proof}

We use Lemma~\ref{lem-prop1} to establish a decay property that is
key to the convergence analysis.
\smallskip

\begin{lemma}\label{L-key-lemma}
Let $(\m{x}^*, \g{\lambda}^*) \in \C{W}^*$ be any solution/multiplier
pair for $(\ref{Prob})$--$(\ref{ProbM})$, let
$\m{x}^k$, $\m{y}^k$, $\m{z}^k$, and $\g{\lambda}^k$ be the iterates
of the generalized BOSVS algorithm, and define
\[
E_k = \rho \| \m{y}_+^k - \m{x}_+^*\|_{\m{P}}^2 +
\frac{1}{\rho} \| \g{\lambda}^k - \g{\lambda}^*\|^2 +  \alpha 
\sum_{i=1}^m \delta_i^k \|\m{x}_i^k - \m{x}_i^*\|^2,
\]
where $\m{P} = \m{MH}^{-1}\m{M}\tr$.
Then for $k$ large enough that the monotonicity condition
$(\ref{delta_monotone})$ holds for all $i \in [1, m]$, we have
\[
E_k \ge E_{k+1} + c_1 \|\m{x}^{k+1}- \m{x}^k\|^2 + c_2\rho
(\|\m{y}_+^k - \m{z}_+^{k}\|_{\m{H}}^2
+ \|\m{Az}^{k} - \m{b}\|^2) ,
\]
where $c_1 = \sigma\alpha\delta_{\min}$ and $c_2 = \alpha(1-\alpha)$.
\end{lemma}
\smallskip

\begin{proof}
By the inequality (\ref{prop1}) of Lemma~\ref{lem-prop1} with
$\m{v} = \m{x}_i^k$, $\m{w} = \m{x}_i^*$, and $\m{u} = \m{z}_i^{k}$, we have
\begin{eqnarray}
L_i^k (\m{x}_i^*)  - L_i^k (\m{z}_i^{k})
- \frac{\rho}{2} \|\m{A}_i\m{z}_{e,i}^{k}\|^2 &\ge&
\frac{\delta_i^k}{2} (\|\m{z}_{e,i}^{k}\|^2
- \|\m{x}_{e,i}^k\|^2)  + 
\frac{\sigma \delta_i^k}{2}\|\m{z}_i^{k}- \m{x}_i^k\|^2 \nonumber \\
&=& \frac{\delta_i^k}{2} (\|\m{x}_{e,i}^{k+1}\|^2
- \|\m{x}_{e,i}^k\|^2)  + 
\frac{\sigma \delta_i^k}{2}\|\m{x}_i^{k+1}- \m{x}_i^k\|^2
\label{L-3456}
\end{eqnarray}
where $\m{x}_e^k = \m{x}^k - \m{x}^*$, $\m{z}_e^k = \m{z}^k - \m{x}^*$,
and $\m{z}^k = \m{x}^{k+1}$ by Step~1b of generalized BOSVS.
Since $\m{x}^*$ minimizes $L(\;\cdot\;, \g{\lambda}^*)$ and since
the augmented Lagrangian is the sum of smooth and nonsmooth terms,
the first-order optimality condition implies that for each $i \in [1,m]$,
\begin{equation}\label{zzz}
\left\langle \m{g}_i^*, \; \m{w} - \m{x}_i^* \right\rangle
+ h_i(\m{w}) - h_i (\m{x}_i^*) \ge 0
\end{equation}
for all $\m{w} \in \mathbb{R}^{n_i}$, where $\m{g}_i^*$ is the gradient
of the smooth part of the objective evaluated at $\m{x}_i^*$:
\[
\m{g}_i^* = \nabla f_i (\m{x}_i^*)  + \rho \m{A}_i  \tr
\left( \sum_{i=1}^m \m{A}_i \m{x}_i^* - \m{b}+ \g{\lambda}^k/\rho \right) .
\]
We add 
$L_i^k (\m{x}_i^*)  - L_i^k (\m{z}_i^{k})
- \rho \|\m{A}_i\m{z}_{e,i}^{k}\|^2/2$ to both sides of
(\ref{zzz}) and take $\m{w} = \m{z}_i^{k}$.
After much cancellation, we obtain the relation
\newpage
\begin{eqnarray}
&L_i^k (\m{x}_i^*)  - L_i^k (\m{z}_i^{k})
- \frac{\rho}{2} \|\m{A}_i\m{z}_{e,i}^{k}\|^2 \le \nonumber \\
&f_i(\m{x}_i^*) - f_i(\m{z}_i^{k}) + \nabla f_i(\m{x}_i^*)\tr
\m{z}_{e,i}^{k} - \rho
\left \langle
\sum_{j \le i} \m{A}_j \m{z}_{e,j}^{k} +
\sum_{j > i} \m{A}_j \m{y}_{e,j}^k 
+ \g{\lambda}_e^k/\rho, \; \m{A}_i \m{z}_{e,i}^{k} \right\rangle \le&
\nonumber \\
&- \rho \left\langle \sum_{j \le i} \m{A}_j \m{z}_{e,j}^{k} +
\sum_{j > i} \m{A}_j \m{y}_{e,j}^k 
+ \g{\lambda}_e^k/\rho, \; \m{A}_i \m{z}_{e,i}^{k} \right\rangle , &
\label{zzz2}
\end{eqnarray}
where $\m{y}_{e}^k = \m{y}^k - \m{x}^*$,
$\g{\lambda}_e^k = \g{\lambda}^k - \g{\lambda}^*$, and
the last inequality is due to the convexity of $f_i$.
We combine this upper bound with the lower bound (\ref{L-3456}) to obtain
\begin{eqnarray} \label{L-key}
&& -\rho \left\langle \m{A}_i \m{z}_{e,i}^{k},
\sum_{j \le i} \m{A}_j \m{z}_{e,j}^{k} + \sum_{j >i} 
\m{A}_j \m{y}_{e,j}^k + \g{\lambda}_e^k/\rho \right\rangle \nonumber \\
&\ge& \frac{\delta_i^k}{2} (\|\m{x}_{e,i}^{k+1}\|^2 - \|\m{x}_{e,i}^k\|^2)
+  \frac{\sigma \delta_i^k}{2}  \| \m{x}_i^{k+1} - \m{x}_i^k\|^2 .
\end{eqnarray}

Focusing on the left side of (\ref{L-key}), observe that
\begin{eqnarray}
\sum_{j \le i} \m{A}_j \m{z}_{e,j}^{k} +
\sum_{j >i} \m{A}_j \m{y}_{e,j}^k &=&
\sum_{j =1 }^m \m{A}_j (\m{z}_{j}^{k} - \m{x}_j^*) + 
\sum_{j >i} \m{A}_j (\m{y}_{j}^k - \m{z}_j^{k}) \nonumber \\
&=& \m{A}\m{z}^{k}-\m{b} + 
\sum_{j >i} \m{A}_j (\m{y}_{j}^k - \m{z}_j^{k}) \label{xyz}
\end{eqnarray}
since $\m{Ax}^* = \m{b}$.
Let $\tau_i^k$ denote the right side of (\ref{L-key}):
\[
\tau_i^k =
\frac{\delta_i^k}{2} (\|\m{x}_{e,i}^{k+1}\|^2 - \|\m{x}_{e,i}^k\|^2)
+  \frac{\sigma \delta_i^k}{2}  \| \m{x}_i^{k+1} - \m{x}_i^k\|^2 .
\]
With this notation and with the simplification (\ref{xyz}),
(\ref{L-key}) becomes
\begin{equation}\label{xzz}
-\rho \left\langle \m{A}_i \m{z}_{e,i}^{k},
\m{A}\m{z}^{k}-\m{b} + \g{\lambda}_e^k/\rho +
\sum_{j >i} \m{A}_j (\m{y}_{j}^k - \m{z}_j^{k}) \right\rangle \ge \tau_i^k .
\end{equation}

We will sum the inequality (\ref{xzz}) over $i$ between 1 and $m$.
Since
\[
\sum_{i = 1}^m \m{A}_i \m{z}_{e,i}^{k} = 
\sum_{i=1}^m \m{A}_i (\m{z}_i^k - \m{x}_i^*) = \m{Az}^k - \m{b} := \m{r}^k ,
\]
it follows that in (\ref{xzz}),
\begin{equation}\label{1st}
\sum_{i = 1}^m \left\langle \m{A}_i \m{z}_{e,i}^{k}, \m{r}^k
+ \g{\lambda}_e^k/\rho \right\rangle =
\left\langle \m{r}^k, \m{r}^k + \g{\lambda}_e^k/\rho \right\rangle.
\end{equation}
Also, observe that
\[
\sum_{j>i} \m{A}_j (\m{y}_{j}^k - \m{z}_j^{k}) =
\sum_{j=2}^m \m{A}_j (\m{y}_{j}^k - \m{z}_j^{k}) -
\sum_{j =2}^i \m{A}_j (\m{y}_{j}^k - \m{z}_j^{k}),
\]
with the convention that the sum from $j =2$ to $j = 1$ is 0.
Take the inner product of this identity with $\m{A}_i\m{z}_{e,i}^{k}$
and sum over $i$ to obtain
\begin{eqnarray}
&\sum_{i=1}^m \left\langle \m{A}_i \m{z}_{e,i}^{k},
\sum_{j >i} \m{A}_j (\m{y}_{j}^k - \m{z}_j^{k}) \right\rangle = & \nonumber \\
&\left\langle \m{r}^k, 
\sum_{j=2}^m \m{A}_j (\m{y}_{j}^k - \m{z}_j^{k}) \right\rangle
- (\m{z}_+^{k} - \m{x}_+^*)\tr \m{M} (\m{y}_+^k -\m{z}_+^{k}) ,& \label{2nd}
\end{eqnarray}
where $\m{M}$ is defined in (\ref{m-def}).
We sum (\ref{xzz}) over $i$ between 1 and $m$ and
utilize (\ref{1st}) and (\ref{2nd}) to obtain
\begin{equation}\label{close}
\quad \quad \;
(\m{y}_+^k - \m{x}_+^*)\tr\m{Mw}
- \frac{1}{\rho} \left( \langle \m{r}^k , \g{\lambda}_e^k \rangle
+ \sum_{i=1}^m \tau_i^k \right) \ge
\m{w}\tr\m{Mw} +
\left\langle
\m{r}^{k} ,
\m{r}^{k} + \sum_{j=2}^m \m{A}_j \m{w}_{j-1}
\right\rangle ,
\end{equation}
where $\m{w} = \m{y}_+^{k} - \m{z}_+^{k}$.

Observe that
\begin{eqnarray*}
\m{w}\tr\m{Mw} &=& \frac{1}{2} \m{w}\tr(\m{M} + \m{M}\tr)\m{w} =
\frac{1}{2} \m{w}\tr(\m{M} + \m{M}\tr - \m{H})\m{w} +
\frac{1}{2}\m{w}\tr\m{Hw} \\
&=& \frac{1}{2}\left\| \sum_{i=2}^m \m{A}_i \m{w}_{i-1}\right\|^2
+ \frac{1}{2} \m{w}\tr\m{Hw}
\end{eqnarray*}
since $(\m{M}+\m{M}\tr - \m{H})_{ij} = \m{A}_{i+1}\tr\m{A}_{j+1}$
by the definition of $\m{M}$ and $\m{H}$.
With this substitution, the right side of (\ref{close})
becomes a sum of squares:
\begin{eqnarray*}
&\m{w}\tr\m{Mw} +
\left\langle
\m{r}^{k} ,
\m{r}^{k} + \sum_{j=2}^m \m{A}_j \m{w}_{j-1}
\right\rangle = &
\frac{1}{2} \left(
\m{w}\tr\m{Hw} + \|\m{r}^k\|^2 +
\left\| \m{r}^k + \sum_{i=2}^m \m{A}_i \m{w}_{i-1}\right\|^2 \right).
\end{eqnarray*}
Hence, it follows from (\ref{close}) that
\begin{equation}\label{closer}
(\m{y}_+^k - \m{x}_+^*)\tr\m{Mw}
- \frac{1}{\rho} \left( \langle \m{r}^k , \g{\lambda}_e^k \rangle
+ \sum_{i=1}^m \tau_i^k \right) \ge
\frac{1}{2} \left( \|\m{w}\|_{\m{H}}^2 + \|\m{r}^k\|^2 \right).
\end{equation}

Let $\m{P} = \m{MH}^{-1} \m{M} \tr$ and recall that $\m{w} =$
$\m{y}_+^{k} - \m{z}_+^{k}$.
By the definition of $\m{y}^{k+1}$ and $\g{\lambda}^{k+1}$ in Step~3
of Algorithm~\ref{ADMMcommon}, we have
\begin{eqnarray*}
& \| \m{y}_+^k - \m{x}_+^*\|_{\m{P}}^2 -
\| \m{y}_+^{k+1}- \m{x}_+^*\|_{\m{P}}^2 +
\frac{1}{\rho^2} (\| \g{\lambda}_e^k\|^2
- \| \g{\lambda}_e^{k+1}\|^2) = & \\
& \| \m{y}_+^k - \m{x}_+^*\|_{\m{P}}^2 -
\| (\m{y}_+^k - \m{x}_+^*) - \alpha  \m{M}^{-\sf T} \m{H}
\m{w} \|_{\m{P}}^2  +
\frac{1}{\rho^2} (\| \g{\lambda}_e^k\|^2
-   \| \g{\lambda}_e^k + \alpha \rho\m{r}^{k}\|^2) = & \\
& 2 \alpha (\m{y}_+^k - \m{x}_+^*)\tr\m{Mw}
- \alpha^2  \|\m{w}\|_{\m{H}}^2 -
\frac{2 \alpha}{\rho} \langle \m{r}^k , \g{\lambda}_e^k \rangle
- \alpha^2 \|\m{r}^k\|^2. &
\end{eqnarray*}
On the right side of this inequality, we utilize
(\ref{closer}) multiplied by $2\alpha$ to conclude that
\begin{eqnarray}\label{drfc}
& \| \m{y}_+^k - \m{x}_+^*\|_{\m{P}}^2 -
\| \m{y}_+^{k+1}- \m{x}_+^*\|_{\m{P}}^2 +
\frac{1}{\rho^2} (\| \g{\lambda}_e^k\|^2
- \| \g{\lambda}_e^{k+1}\|^2) - \frac{2\alpha}{\rho}
\sum_{i=1}^m \tau_i^k \ge
\nonumber \\
& c_2 (\|\m{y}_+^{k} - \m{z}_+^{k}\|_{\m{H}}^2 + \|\m{r}^k\|^2)
\end{eqnarray}
where $c_2 = \alpha (1-\alpha)  > 0$ since $\alpha \in (0,1)$.
By the definition of $\tau_i^k$ and the assumption that
$k$ is large enough that (\ref{delta_monotone}) holds for all $i$,
it follows that
\begin{eqnarray*}
-\tau_i^k &=&
\frac{\delta_i^k}{2} (\|\m{x}_{e,i}^{k}\|^2 - \|\m{x}_{e,i}^{k+1}\|^2)
-  \frac{\sigma \delta_i^k}{2}  \| \m{x}_i^{k+1} - \m{x}_i^k\|^2 \\
&\le&
\frac{\delta_i^k}{2} \|\m{x}_{e,i}^{k}\|^2 -
\frac{\delta_i^{k+1}}{2} \|\m{x}_{e,i}^{k+1}\|^2
-  \frac{\sigma \delta_i^k}{2}  \| \m{x}_i^{k+1} - \m{x}_i^k\|^2 .
\end{eqnarray*}
This bound for $-\tau_i^k$ along with the inequality (\ref{drfc}) and
the definition of $E^k$ complete the proof.
\end{proof}

The following theorem establishes the global convergence of generalized BOSVS.
\smallskip

\begin{theorem}\label{L-glob-thm}
If $\m{x}^k$, $\m{y}^k$, and $\g{\lambda}^k$ are iterates of the generalized
BOSVS algorithm, then the $\m{x}^k$ and
$\m{y}^k$ sequences converge to a common limit denoted $\m{x}^*$ and
the $\g{\lambda}^k$ converge to a limit denoted $\g{\lambda}^*$ where
$(\m{x}^*, \g{\lambda}^*) \in \C{W}^*$.
\end{theorem}
\smallskip

\begin{proof}
Let $\bar{k}$ be chosen large enough that (\ref{delta_monotone}) holds
for all $k \ge \bar{k}$.
Since $\m{x}^{k+1} = \m{z}^k$ in generalized BOSVS,
it follows from Lemma~\ref{L-key-lemma} that for
$j \ge \bar{k}$ and $p > 0$, we have
\begin{equation}\label{Ej}
\quad \quad
E_j \ge E_{j+p} + c \sum_{k=j}^{j+p-1}
(\|\m{x}^{k+1}- \m{x}^k\|^2 +
\|\m{y}_+^k - \m{x}_+^{k+1}\|_{\m{H}}^2
+ \|\m{Ax}^{k+1} - \m{b}\|^2 ),
\end{equation}
where $c = \min\{c_1, \rho c_2\} > 0$.
Let $p$ tend to $+\infty$ in (\ref{Ej}).
Since the columns of $\m{A}_i$ are linearly independent for $i \ge 2$,
$\m{H}$ is positive definite and
\begin{equation}\label{L-lim1}
\lim_{k \to \infty} \|\m{x}^{k+1}- \m{x}^k\|=
\lim_{k \to \infty} \|\m{y}_+^{k}- \m{x}_+^{k+1}\|=
\lim_{k \to \infty} \|\m{Ax}^{k+1}- \m{b}\|= 0.
\end{equation}
By the definition of $\m{b}_i^k$, we have
\begin{equation}\label{Ax-blim}
\m{A}_i \m{x}_i^{k+1} -\m{b}_i^k =
\sum_{j \le i }  \m{A}_j \m{x}_j^{k+1} +
\sum_{j > i } \m{A}_j \m{y}_j^k - \m{b}.
\end{equation}
By (\ref{L-lim1}), $\m{y}_+^k$ approach $\m{x}_+^{k+1}$, and by (\ref{Ax-blim})
and (\ref{L-lim1}),
\begin{equation}\label{L-lim-b}
\lim_{k \to \infty} (\m{A}_i \m{x}_i^{k+1} -\m{b}_i^k) =
\lim_{k \to \infty}
\m{A} \m{x}^{k+1} -\m{b} = \m{0}
\end{equation}
for all $i \in [1,m]$.

By the definition of $E_k$ in Lemma~\ref{L-key-lemma}, we see that the
iterates $\g{\lambda}^k$ and $\m{x}^k$ are uniformly bounded.
Hence, there exist limits $\g{\lambda}^*$ and $\m{x}^*$, and an infinite
sequence $\C{K} \subset \{1, 2, \ldots \}$ such that
$\g{\lambda}^k$ and $\m{x}^k$ for $k \in \C{K}$ converge to $\g{\lambda}^*$
and $\m{x}^*$ respectively.
By the first relation in (\ref{L-lim1}),
$\m{x}^{k+1}$ also converges to $\m{x}^*$ for $k \in \C{K}$.
In Step~1b of generalized BOSVS, we have
\[
\m{x}_i^{k+1} = \arg \min \{
\Phi_i^k (\m{u}, \m{x}_i^k, \delta_i^k) : \m{u} \in \mathbb{R}^{n_i} \} .
\]
The first-order optimality conditions for $\m{x}_i^{k+1}$ are
\begin{equation}\label{1storder}
\left\langle \m{g}_i^k, \; \m{u} - \m{x}_i^{k+1} \right\rangle
+ h_i(\m{u}) \ge h_i (\m{x}_i^{k+1})
\end{equation}
for all $\m{u} \in \mathbb{R}^{n_i}$, where $\m{g}_i^k$ is the gradient
of the smooth part of the objective evaluated at $\m{x}_i^{k+1}$:
\[
\m{g}_i^k =
\nabla f_i (\m{x}_i^{k})  + \rho \m{A}_i  \tr (\m{A}_i \m{x}_i^{k+1} -
\m{b}_i^k+ \g{\lambda}^k/\rho) + \delta_i^k (\m{x}_i^{k+1} - \m{x}_i^{k}) .
\]
As $k \in \C{K}$ tends to infinity, $\nabla f_i (\m{x}_i^{k})$ approaches
$\nabla f_i(\m{x}^*)$ since $\nabla f_i$ is Lipschitz continuous,
$\m{A}_i \m{x}_i^{k+1} - \m{b}_i^k$ approaches $\m{0}$ by (\ref{L-lim-b}), and
$\delta_i^k (\m{x}_i^{k+1} - \m{x}_i^{k})$ approaches $\m{0}$ by (\ref{L-lim1})
and the uniform bounded (\ref{delta_monotone}) for $\delta_i^k$.
Consequently, we have
\begin{equation}\label{glim}
\lim_{k \in \C{K}} \; \m{g}_i^k =
\nabla f_i(\m{x}^*) + \m{A}_i\tr \g{\lambda}^* .
\end{equation}
Let $k \in \C{K}$ tend to $+\infty$ in (\ref{1storder}).
By (\ref{glim}) and the lower semicontinuity of $h_i$, we deduce that
\[
\langle \nabla f_i (\m{x}_i^*) + \m{A}_i \tr \g{\lambda}^*,
\m{u} - \m{x}_i^* \rangle + h_i(\m{u}) \ge h_i (\m{x}_i^*)
\]
for all $\m{u} \in \mathbb{R}^{n_i}$.
Therefore, $\m{x}^*$ and $\g{\lambda}^*$ satisfy the first-order
optimality condition (\ref{Wstar}).
By the last relation in (\ref{L-lim1}), it follows that $\m{Ax}^* = \m{b}$
and $\m{x}^*$ is feasible in (\ref{Prob}).
By the convexity of $f_i$ and $h_i$,
$\m{x}^*$ is a solution of (\ref{Prob})--(\ref{ProbM})
and $\g{\lambda}^*$ is an associated multiplier for the linear constraint.

Since $\m{x}^{k+1}$ converges to $\m{x}^*$ for $k \in \C{K}$,
the second relation in (\ref{L-lim1}) implies that
$\m{y}_+^k$ converges to $\m{x}_+^*$ for $k \in \C{K}$.
In Lemma~\ref{L-key-lemma}, we use the specific limits
$\m{x}^*$ and $\g{\lambda}^*$ associated with $k \in \C{K}$.
Hence, $E_k$ tends to 0 for $k \in \C{K}$.
It follows from (\ref{Ej}) that the entire $E_k$ sequence tends to 0.
By the definition of $E_k$, we deduce that
the entire $(\m{x}^k, \m{y}_+^k, \g{\lambda}^k)$ sequence converges
$(\m{x}^*, \m{x}_+^*, \g{\lambda}^*)$.
Since $\m{y}_1^k = \m{x}_1^k$ for each $k$, where $\m{x}_1^k$ converges
to $\m{x}_1^*$, we conclude that $\m{y}^k$ converges to $\m{x}^*$.
This completes the proof.
\end{proof}
\section{Multistep BOSVS}
\label{mBOSVS}
For the template given by Algorithm~\ref{ADMMcommon},
we only need to assume that the columns of $\m{A}_i$ are linearly
independent for $i \ge 2$ since only these columns enter into the matrix
$\m{M}$ which is inverted in Step~3.
For generalized BOSVS, this assumption was sufficient of convergence.
On the other hand, for both multistep and accelerated BOSVS,
strong convexity of the augmented Lagrangian with respect to each
of the variables $\m{x}_i$ is needed in the analysis.
Since it has already been assumed that the columns of $\m{A}_i$ are
linearly independent for $i \ge 2$, we will simply strengthen this
assumption to require, henceforth, that the columns
of $\m{A}_i$ are linearly independent for every $i$.
This ensures strong convexity of the augmented Lagrangian $L$ with
respect to each of the variables $\m{x}_i$.

The inner loop for the multistep BOSVS algorithm appears in Algorithm~\ref{2}.
\renewcommand\figurename{Alg.}
\begin{figure}[h]
{\tt
\begin{tabular}{l}
\hline
{\bf Initialize:} $\m{u}_i^0 = \m{x}_i^k$ \\[.05in]
{\bf For } $l = 1, 2, \ldots $\\[.05in]
\begin{tabular}{ll}
1a. & Choose $\delta_0^l \in [\delta_{\min}, \delta_{\max}]$. \\[.05in]
1b. & Set $\delta^l = \eta^j \delta_{0}^l$, where
$j \ge 0$ is the smallest integer such that \\[.05in]
& $\quad$ $f_i(\m{u}_i^{l-1}) +$
$\langle \nabla f_i(\m{u}_i^{l-1}), \m{u}_i^{l} -\m{u}_i^{l-1} \rangle
+ \frac{(1-\sigma)\delta^l}{2} \|\m{u}_i^{l} - \m{u}_i^{l-1}\|^2
\ge f_i(\m{u}_i^{l})$,\\[.05in]
&
where
$\m{u}_i^{l} = \arg \min \{
\Phi_i^k (\m{u}, \m{u}_i^{l-1}, \delta_i^k) : \m{u} \in \mathbb{R}^{n_i} \}$.
\\[.05in]
1c. &
If $\gamma^l := \sum_{j=1}^l 1/\delta^j \ge \Gamma_i^{k-1}$ and 
$\| \m{u}_i^l - \m{u}_i^{l-1} \|/\sqrt{\gamma^l} \le \psi (e^{k-1})$,
break.\\[.05in]
\end{tabular}\\
{\bf Next} \\[.05in]
1d. Set
$\m{z}_i^k = \left( \sum_{j=1}^l \m{u}_i^j/\delta^j \right)/\gamma^l$,
$r_i^k = (1/\gamma^l) \sum_{j=1}^l \|\m{u}_i^j - \m{u}_i^{j-1}\|^2$,
\\[.05in]
$\quad \quad \quad \m{x}_i^{k+1}= \m{u}_i^l$, and $\Gamma_i^k = \gamma^l$.\\
\hline
\end{tabular}
}
\caption{Inner loop in Step~$1$ of Algorithm~$\ref{ADMMcommon}$
for the multistep BOSVS scheme.}
\label{2} 
\end{figure}
\renewcommand\figurename{Fig.}
In generalized BOSVS, the iteration is given by
$\m{x}_i^{k+1} =$
$\arg \min \{ \Phi_i^k (\m{u}, \m{x}_i^k, \delta_i^k) :
\m{u} \in \mathbb{R}^{n_i} \}$ where $\delta_i^k$ is
determined by a line search process.
In the multistep BOSVS algorithm, this
single minimization is replaced by the recurrence
\[
\m{u}_i^{l} =
\arg \min \{ \Phi_i^k (\m{u}, \m{u}_i^{l-1}, \delta_i^k) :
\m{u} \in \mathbb{R}^{n_i} \},
\]
where $\m{u}_i^0 = \m{x}_i^k$.
By converting the single minimization into a recurrence,
we hope to a achieve a better minimizer of the augmented Lagrangian.
In generalized BOSVS, the convergence relies on a careful choice
of $\delta_i^k$ using both the BB-formula and a safeguarding technique.
In multistep BOSVS, these restrictions on $\delta_i^k$ are replaced
in Step~1c by a condition related to the accuracy of the iterates.
In this step, $\psi: \mathbb{R}^+ \to \mathbb{R}^+$ denotes any function
satisfying $\psi (0) = 0$ and $\psi (s) > 0$ for $s > 0$ with $\psi$
continuous at $s = 0$.
For example, $\psi(t) = t$.

\begin{remark}
Computationally, it is not necessary to store the
$\m{u}_i^l$ sequence to evaluate $\m{z}_i^k$ in Step~{\rm 1d}.
For example, in Step~{\rm 1c} we could introduce a sequence
\[
\m{a}_i^l = (1-\alpha^l) \m{a}_i^{l-1} + \alpha^l \m{u}_i^l, \quad
\mbox{\rm where } \alpha^l = 1/(\delta^l \gamma^l)
\mbox{ \rm and } \m{a}_i^0 = \m{x}_i^0,
\]
and in Step~{\rm 1d}, we would set $\m{z}_i^k = \m{a}_i^l$.
Note that $0 < \alpha^l \le 1$ due to the form of
$\gamma^l$ in Step~{\rm 1c}.
\end{remark}

Since $\eta > 1$, the line search in Step~1b of multistep BOSVS
terminates in a finite number of iterations and the final $\delta^l$
has exactly the same bounds (\ref{delta_bound}) as that of
generalized BOSVS.
Since $\delta^l$ is uniformly bounded, it follows that the condition
$\Gamma_i^k \ge \Gamma_i^{k-1}$ of Step~1c is fulfilled for $l$ sufficiently
large.
In the numerical experiments for multistep BOSVS in Section~\ref{numerical},
$\delta_0^l$ is given by the safeguarded BB choice of generalized BOSVS.

Let us first observe that when
$e^k =0$, we have reached a solution of (\ref{Prob})--(\ref{ProbM}).
\smallskip

\begin{lemma}\label{L-stop-cond2}
If $e^{k}=0$ in the multistep BOSVS algorithm,
then $\m{x}^{k+1} = \m{x}^k = \m{y}^k = \m{z}^k$ solves
$(\ref{Prob})$--$(\ref{ProbM})$ and $(\m{x}^k, \g{\lambda}^k) \in \C{W}^*$.
\end{lemma}
\smallskip

\begin{proof}
If $e^k = 0$, then $r_i^k = 0$ for each $i$.
It follows that $\m{x}_i^k = \m{u}_i^0 = \m{u}_i^1 = \ldots = \m{u}_i^l$.
By Step~1d, $\m{u}_i^l = \m{x}_i^{k+1} = \m{z}_i^k$.
Consequently, we have $\m{x}^{k+1} = \m{x}^k = \m{z}^k$.
Since all three algorithms in this paper share
Algorithm~\ref{ADMMcommon}, the remainder of the proof is
exactly as in Lemma~\ref{L-stop-cond}.
\end{proof}

The following inequality is based on Lemma~\ref{lem-prop1}.
\begin{lemma}\label{lem-GD-Conv}
In multistep BOSVS, we have
\begin{equation}\label{GD-Converge}
\nu_i \rho \| \m{z}_i^{k} - \bar{\m{x}}_i^k\|^2 
+ \frac{\sigma}{\Gamma_i^{k}}
\sum_{l=1}^{l_i^k}  \|\m{u}_i^l- \m{u}_i^{l-1}\|^2 
  \le  \frac{\| \m{x}_i^k - \bar{\m{x}}_i^k\|^2}
{\Gamma_i^{k}} ,
\end{equation}
for each $i \in [1,m]$,
where $l_i^k$ is the terminating value of $l$ at iteration $k$,
$\nu_i>0$ is the smallest eigenvalue of $\m{A}_i \tr \m{A}_i$, and
\begin{equation}\label{xbar}
\bar{\m{x}}_i^k= \arg \min \{ L_i^k(\m{u}) : \m{u} \in \mathbb{R}^{n_i}\}
\end{equation}
with $L_i^k$ defined in $(\ref{lik})$.
\end{lemma}
\smallskip

\begin{proof}
By Lemma \ref{lem-prop1}, we have
\begin{eqnarray} \label{GD-linesearch}
L_i^k (\m{w})  - L_i^k (\m{u}_i^l) &\ge&
\frac{\delta^l}{2} (\|\m{w} - \m{u}_i^l\|^2 - 
\|\m{w} - \m{u}_i^{l-1}\|^2 )
+ \frac{\rho}{2} \| \m{A}_i (\m{w} - \m{u}_i^l)\|^2  \\
&& + \frac{\sigma\delta^l}{2}  \|\m{u}_i^l- \m{u}_i^{l-1}\|^2
\nonumber
\end{eqnarray}
for any $\m{w} \in \mathbb{R}^{n_i}$.
We take $\m{w} = \bar{\m{x}}_i^k$.
Since $L_i^k (\m{u}_i^l) - L_i^k (\bar{\m{x}}_i^k) \ge 0$,
we have
\begin{equation}\label{xbarbound}
\frac{\rho}{\delta^l} \| \m{A}_i (\bar{\m{x}}_i^k - \m{u}_i^l)\|^2
+ \sigma \|\m{u}_i^l- \m{u}_i^{l-1}\|^2
\le \|\bar{\m{x}}_i^k - \m{u}_i^{l-1}\|^2 - \|\bar{\m{x}}_i^k
- \m{u}_i^l\|^2 .
\end{equation}
Summing this inequality for $l$ between 1 and $l_i^k$ gives
\begin{equation}\label{tyhg}
\rho \sum_{l=1}^{l_i^k} \frac{1}{\delta^l}
\| \m{A}_i (\bar{\m{x}}_i^k - \m{u}_i^l)\|^2 
+ \sigma \sum_{l=1}^{l_i^k}  \|\m{u}_i^l- \m{u}_i^{l-1}\|^2 \le
\|\bar{\m{x}}_i^k - \m{x}_i^k\|^2 .
\end{equation}
Since the quadratic $\|\m{A}_i (\bar{\m{x}}_i^k - \m{u})\|^2$ is
a convex function of $\m{u}$, it follows from Jensen's inequality that
\[
\sum_{l=1}^{l_i^k} \frac{1}{\delta^l} \| \m{A}_i (\bar{\m{x}}_i^k
- \m{u}_i^l)\|^2 \ge 
\Gamma_i^k   \| \m{A}_i (\bar{\m{x}}_i^k - \m{z}_i^k)\|^2 \ge 
\Gamma_i^k \nu_i  \| \m{z}_i^k - \bar{\m{x}}_i^k\|^2,
\]
where $\nu_i>0$ is the smallest eigenvalue of $\m{A}_i \tr \m{A}_i$.
Combine this with (\ref{tyhg}) to obtain (\ref{GD-Converge}).
\end{proof}

\begin{remark}\label{remark1}
$L_i^k$ is strongly convex since it is the sum of convex functions
and a strongly convex quadratic $\langle \m{A}_i\m{u}, \m{A}_i \m{u}\rangle$;
consequently, the minimizer $\bar{\m{x}}_i^k$ exists.
Due to the upper bound {\rm (\ref{delta_bound})}
for $\delta^l$ in multistep BOSVS,
$\gamma^l$ grows linearly in $l$.
Hence, for the inner loop of multistep BOSVS,
{\rm (\ref{GD-Converge})} implies that
$\| \m{z}_i^{k} - \bar{\m{x}}_i^k\| = O(1/\sqrt{l_i^k})$.
By $(\ref{GD-linesearch})$, the objective values satisfy
$L_i^k(\m{z}_i^k) - L_i^k(\bar{\m{x}}_i^k) = O(1/l_i^k)$;
to see this, divide $(\ref{GD-linesearch})$ by $\delta^l$,
sum over $l$ between 1 and $l_i^k$, and apply Jensen's
inequality twice, to the terms involving $L(\m{u}^l)$ and to the
terms involving $\m{A}_i$.
\end{remark}

As a consequence of Lemma~\ref{lem-GD-Conv}, we show that the stopping
conditions in Step~1c of multistep BOSVS are satisfied for a finite $l$.

\begin{corollary}\label{finite_stop}
If $e^{k-1} >0$ in Step~$1${\rm d} of multistep BOSVS,
then the stopping condition in Step~$1${\rm c} is fulfilled when $l$ is
sufficiently large.
\end{corollary}
\smallskip

\begin{proof}
Since $\delta^l$ in multistep BOSVS has the same upper bound
(\ref{delta_bound}) as generalized BOSVS, it follows that $\Gamma_i^k$
in step~1c of multistep BOSVS tends to infinity as $l$ tends to infinity.
By Lemma~\ref{lem-GD-Conv}, the iteration difference
$\|\m{u}_i^l - \m{u}_i^{l-1}\|$ tends to zero as $l_i^k$ grows.
Hence, both conditions in Step~1c are satisfied for $l$ sufficiently large
when $\psi (e^{k-1}) > 0$.
\end{proof}

Similar to generalized BOSVS, the key to the convergence of multistep
BOSVS is a decay property for the iterates.
The analogue of Lemma~\ref{L-key-lemma} for multistep BOSVS is the
following result.
\begin{lemma}\label{L-key-lemma2}
Let $(\m{x}^*, \g{\lambda}^*) \in \C{W}^*$ be any solution/multiplier
pair for $(\ref{Prob})$--$(\ref{ProbM})$, let
$\m{x}^k$, $\m{y}^k$, $\m{z}^k$, $\m{u}_{k}^l$, and $\g{\lambda}^k$
be the iterates of the multistep BOSVS algorithm,
let $l_i^k$ be the terminating value of $l$ at iteration $k$, and define
\[
E_k = \rho \| \m{y}_+^k - \m{x}_+^*\|_{\m{P}}^2 +
\frac{1}{\rho} \| \g{\lambda}^k - \g{\lambda}^*\|^2 +  \alpha 
\sum_{i=1}^m  \frac{\|\m{x}_i^k - \m{x}_i^*\|^2}{\Gamma_i^k},
\]
where $\m{P} = \m{MH}^{-1}\m{M}\tr$.
Then for all $k$ and for all $i \in [1, m]$, we have
\[
E_k \ge E_{k+1} + c_1 \sum_{i=1}^m \sum_{l=1}^{l_i^k}
\frac{\|\m{u}_{i,k}^l- \m{u}_{i,k}^{l-1}\|^2}{\Gamma_i^k} + c_2\rho
(\|\m{y}_+^k - \m{z}_+^{k}\|_{\m{H}}^2
+ \|\m{Az}^{k} - \m{b}\|^2),
\]
where $c_1 = \sigma\alpha$ and $c_2 = \alpha(1-\alpha)$.
\end{lemma}
\smallskip

\begin{proof}
We put $\m{w} = \m{x}_i^*$ in (\ref{GD-linesearch}) to obtain
\[
\frac{ L_i^k (\m{x}_i^*) - F_i^k (\m{u}_{i,k}^l)}
{\delta^l} \ge 
\frac{1}{2} (\|\m{x}_i^* - \m{u}_{i,k}^l\|^2
- \|\m{x}_i^* - \m{u}_{i,k}^{l-1}\|^2)  + 
\frac{\sigma}{2} \|\m{u}_{i,k}^l- \m{u}_{i,k}^{l-1}\|^2,
\]
where $F_i^k(\m{u}_{i,k}^l) =$
$L_i^k(\m{u}_{i,k}^l) + (\rho/2) \|\m{A}_i( \m{u}_{i,k}^l - \m{x}_i^*) \|^2$.
Summing this inequality over $l$ yields
\begin{eqnarray} \label{2wsx-1}
&& \sum_{l=1}^{l_i^k}
\left( \frac{ L_i^k (\m{x}_i^*) - F_i^k (\m{u}_{i,k}^l)}
{\delta^l} \right) \ge \nonumber \\
&& \frac{1}{2} (\|\m{x}_i^* - \m{u}_{i,k}^{l_i^k}\|^2 - \|\m{x}_i^* -
\m{u}_{i,k}^0\|^2)  + 
\frac{\sigma}{2} \sum_{l=1}^{l_i^k} \|\m{u}_{i,k}^l-
\m{u}_{i,k}^{l-1}\|^2.
\end{eqnarray}
Since $F_i^k$ is convex, it follows from Jensen's inequality and the
definition of $\Gamma_i^k$ and $\m{z}_i^k$ in Step~1c of multistep BOSVS that
\begin{equation}\label{jensen}
\frac{1}{\Gamma_i^k} \sum_{l=1}^{l_i^k}
\frac{1}{\delta^l} F_i^k (\m{u}_{i,k}^l ) \ge
F_i^k \left(\frac{1}{\Gamma_i^k}
\sum_{l=1}^{l_i^k} \frac{1}{\delta^l} \m{u}_{i,k}^l \right)
=  F_i^k (\m{z}_i^k).
\end{equation}
Substitute $\m{x}_i^{k+1} = \m{u}_{i,k}^{l_i^k}$ and
$\m{x}_i^k = \m{u}_{i,k}^0$ in (\ref{2wsx-1}) and use
(\ref{jensen}) to obtain
\begin{equation}\label{jumppoint}
L_i^k (\m{x}_i^*)  - F_i^k (\m{z}_i^k)
\ge \frac{1}{2\Gamma_i^k} (\|\m{x}_{e,i}^{k+1}\|^2 - \|\m{x}_{e,i}^k\|^2)
+ \frac{\sigma}{2\Gamma_i^k}
\sum_{l = 1}^{l_i^k} \|\m{u}_{i,k}^{l} - \m{u}_{i,k}^{l-1}\|^2 ,
\end{equation}
where $\m{x}_{e,i}^k = \m{x}_i^k - \m{x}_i^*$.
By (\ref{zzz2}), we have the upper bound
\[
L_i^k (\m{x}_i^*)  - F_i^k (\m{z}_i^k) \le
- \rho \left\langle \sum_{j \le i} \m{A}_j \m{z}_{e,j}^{k} +
\sum_{j > i} \m{A}_j \m{y}_{e,j}^k 
+ \g{\lambda}_e^k/\rho, \; \m{A}_i \m{z}_{e,i}^{k} \right\rangle .
\]
Combining lower and upper bounds gives
\begin{eqnarray}
&- \rho \left\langle \sum_{j \le i} \m{A}_j \m{z}_{e,j}^{k} +
\sum_{j > i} \m{A}_j \m{y}_{e,j}^k
+ \g{\lambda}_e^k/\rho, \; \m{A}_i \m{z}_{e,i}^{k} \right\rangle \ge &
\label{combine} \\
&\frac{1}{2\Gamma_i^k} (\|\m{x}_{e,i}^{k+1}\|^2 - \|\m{x}_{e,i}^k\|^2)
+ \frac{\sigma}{2\Gamma_i^k}
\sum_{l = 1}^{l_i^k} \|\m{u}_{i,k}^{l} - \m{u}_{i,k}^{l-1}\|^2 , & \nonumber
\end{eqnarray}
which is the same as (\ref{L-key}) but with the following exchanges:
\[
\delta_i^k \longleftrightarrow
1/\Gamma_i^k \quad \mbox{and}\quad
\|\m{x}_i^{k+1} - \m{x}_i^k\|^2 \longleftrightarrow
\sum_{l = 1}^{l_i^k} \|\m{u}_{i,k}^{l} - \m{u}_{i,k}^{l-1}\|^2.
\]
Except for these adjustments,
the remainder of the proof is the same as the proof of
Lemma~\ref{L-key-lemma}, starting with equation (\ref{xyz}).
\end{proof}

Using Lemma~\ref{L-key-lemma2}, we can now prove the convergence of
multistep BOSVS.
The analysis parallels that of Theorem~\ref{L-glob-thm}.
To facilitate the analysis, we recall the definition and some properties
of the prox function.
For any closed convex real-valued function $h$,
\[
\prox_h (\m{v}) = \arg \min \left\{ h(\m{u})
+ \frac{1}{2} \|\m{v} - \m{u}\|^2 : \m{u} \in \mbox{dom}(h) \right\} .
\]
As shown in \cite[p.~340]{Rockafellar70},
the prox function is nonexpansive:
\[
\|\prox_h(\m{v}_1) - \prox_h(\m{v}_2)\| \le \|\m{v}_1 - \m{v}_2\|.
\]
Moreover, if $g$ is a differentiable convex function and
\begin{equation}\label{pre_prox_min}
\m{u}^* = \arg \min_{\m{u}} g(\m{u}) + h(\m{u}),
\end{equation}
then it follows from the first-order optimality conditions for $\m{u}^*$ that
\begin{equation}\label{prox_min}
\m{u}^* = \prox_h (\m{u}^* - \nabla g(\m{u}^*)).
\end{equation}
Conversely, if (\ref{prox_min}) holds, then so does (\ref{pre_prox_min}).
Hence, these relations are equivalent.
These properties will be used in the convergence analysis of multistep BOSVS.
\begin{theorem}\label{L-glob-thm2}
If multistep BOSVS performs an infinite number of iterations
generating iterates $\m{y}^k$, $\m{z}^k$, and $\g{\lambda}^k$,
then the sequences $\m{y}^k$ and $\m{z}^k$ both approach a common
limit $\m{x}^*$ and $\g{\lambda}^k$ approaches a limit $\g{\lambda}^*$ where
$(\m{x}^*, \g{\lambda}^*) \in \C{W}^*$.
\end{theorem}
\smallskip

\begin{proof}
For any $p > 0$, we sum the decay property of Lemma~\ref{L-key-lemma2} 
to obtain
\begin{equation}\label{Ej2}
\quad \quad \,
E_j \ge E_{j+p} + c \sum_{k=j}^{j+p-1} \left(
\|\m{y}_+^k - \m{z}_+^{k}\|_{\m{H}}^2
+ \|\m{Az}^{k} - \m{b}\|^2 +
\sum_{i=1}^m \sum_{l=1}^{l_i^k}
\frac{\|\m{u}_{i,k}^l- \m{u}_{i,k}^{l-1}\|^2}{\Gamma_i^k} \right) , \\
\end{equation}
where $c = \min\{c_1, \rho c_2\} > 0$.
Let $p$ tend to $+\infty$.
Since $\m{H}$ is positive definite, and the $\Gamma_i^k$ are monotone
nondecreasing as a function of $k$,
it follows from (\ref{Ej2}) that
\begin{equation}\label{L-lim12}
\lim_{k \to \infty} \|\m{y}_+^{k}- \m{z}_+^{k}\|= 0 =
\lim_{k \to \infty} \|\m{Az}^{k}- \m{b}\|.
\end{equation}
Moreover, by the definition of $E_k$ in Lemma~\ref{L-key-lemma2},
$\m{y}_+^k$ and $\g{\lambda}^k$ are bounded sequences, and by the
first equation in (\ref{L-lim12}), $\m{z}_+^k$ is also a bounded sequence.
The second equation in (\ref{L-lim12}) is equivalent to
\[
\lim_{k \to \infty} \left\| \m{A}_1 \m{z}_1^k -
\left( \m{b} - \sum_{i=2}^m \m{A}_i \m{z}_i^k \right) \right\| = 0.
\]
Since $\m{z}_+^k$ is bounded and the columns of $\m{A}_1$ are linearly
independent, $\m{z}_1$ is bounded.
Hence, both $\m{z}^k$ and $\g{\lambda}^k$ are bounded sequences, and
there exist an infinite sequence $\C{K} \subset \{1, 2, \ldots \}$
and limits $\m{x}^*$ and $\g{\lambda}^*$ such that
\begin{equation}\label{limz}
\lim_{k\in\C{K}} \m{z}^k = \m{x}^* \quad \mbox{and} \quad
\lim_{k\in\C{K}} \g{\lambda}^k = \g{\lambda}^*.
\end{equation}
By the first equation in (\ref{L-lim12}), we have
\begin{equation}\label{limy}
\lim_{k \in \C{K}} \m{y}_+^k = \m{x}_+^*.
\end{equation}
By the second equation in (\ref{L-lim12}), $\m{Ax}^* = \m{b}$.
Consequently, by (\ref{limz}) and (\ref{limy}),
\begin{equation}\label{L-lim-b2}
\lim_{k \in \C{K}} \left( \m{A}_i \m{z}_i^{k} -\m{b}_i^k \right) =
\lim_{k \in \C{K}} \left(
\sum_{j \le i }  \m{A}_j \m{z}_j^{k} +
\sum_{j > i } \m{A}_j \m{y}_j^k - \m{b} \right) = \m{Ax}^* - \m{b} = \m{0}
\end{equation}
for all $i \in [1,m]$.

The decay property (\ref{Ej2}) also implies that for each $i$,
\begin{equation}\label{riklim}
\lim_{k \to \infty} r_i^k =
\lim_{k \to \infty} \frac{1}{\Gamma_i^k} \sum_{l=1}^{l_i^k}
\|\m{u}_{i,k}^l- \m{u}_{i,k}^{l-1}\|^2 = 0.
\end{equation}
Combine this with (\ref{L-lim12}) to conclude that
\begin{equation}\label{eklim}
\lim_{k \to \infty} e^k = \lim_{k \to \infty} \psi(e^k) = 0.
\end{equation}
%
The remainder of the proof is partitioned into two cases depending on whether
the monotone nondecreasing sequence $\Gamma_i^k$
either approaches a finite limit, or tends to infinity.

{\bf Case 1.} For some $i$, $\Gamma_i^k$ approaches a finite limit.
Due to the upper bound (\ref{delta_bound}) for $\delta^l$ in Step~1b
of multistep BOSVS, we conclude that $l_i^k$ is uniformly bounded.
By (\ref{riklim}), $\| \m{u}_{i,k}^l - \m{u}_{i,k}^{l-1} \|$ approaches zero,
where the convergence is uniform in $k$ and $l \in [1, l_i^k]$.
Since $\m{u}_{i,k}^0 = \m{x}_i^k$, the triangle inequality and the
uniform upper bound for $l_i^k$ imply that
$\| \m{x}_i^k - \m{u}_{i,k}^l \|$ approaches zero,
where the convergence is uniform in $k$ and $l \in [1, l_i^k]$.
Since $\m{z}_i^k$ is a convex
combination of $\m{u}_{i,k}^l$ for $0 \le l \le l_i^k$ with $l_i^k$
uniformly bounded and $\| \m{x}_i^k - \m{u}_{i,k}^l \|$ approaching zero,
it follows that $\|\m{z}_i^k - \m{x}_i^k\|$ approaches zero.
We summarize these observations in the relation
\begin{equation}\label{uiklim}
\lim_{k \to \infty} \|\m{z}_i^k - \m{x}_i^k\| =
\lim_{k \to \infty} \|\m{z}_i^k - \m{u}_{i,k}^0\| =
\lim_{k \to \infty} \|\m{z}_i^k - \m{u}_{i,k}^1\| = 0 .
\end{equation}
%

In multistep BOSVS,
$\m{u}_{i,k}^1$ minimizes $\Phi_i(\cdot, \m{u}_i^0, \delta_i^k)$.
Identify $g$ in (\ref{pre_prox_min}) with the smooth terms in $\Phi_i$.
By (\ref{prox_min}), we have
\[
\m{u}_{i,k}^1 = \prox_{h_i} \left(
\m{u}_{i,k}^1 - \nabla f_i (\m{u}_{i,k}^0) - \delta_i^k (\m{u}_{i,k}^1
- \m{u}_{i,k}^0 ) - \rho \m{A}_i\tr (\m{A}_i \m{u}_{i,k}^1 - \m{b}_i^k
+ \g{\lambda}^k/\rho) \right) .
\]
Let us now take the limit as $k$ tends to infinity with $k \in \C{K}$.
By (\ref{limz}), $\m{z}_i^k$ approaches $\m{x}_i^*$.
By (\ref{uiklim}) both $\m{u}_{i,k}^0$ and
$\m{u}_{i,k}^1$ approach $\m{z}_i^k$,
and by (\ref{L-lim-b2}) $\m{A}_i \m{u}_{i,k}^1 - \m{b}_i^k$ approaches zero.
Since the prox function and $\nabla f_i$ are both Lipschitz continuous,
we deduce that in the limit, as $k$ tends to infinity with $k \in \C{K}$,
\[
\m{x}_i^* = \prox_{h_i} \left( \m{x}_i^* - \nabla f_i
(\m{x}_{i}^*) - \m{A}_i\tr \g{\lambda}^* \right) .
\]
Again, by (\ref{pre_prox_min}), we have
\begin{equation}\label{stationary1}
\m{x}_i^* = \arg \min \{ f_i(\m{u}) + h_i (\m{u}) +
\langle \g{\lambda}^* , \m{A}_i \m{u} \rangle : \m{u} \in \mathbb{R}^{n_i} \}.
\end{equation}
If this were to hold for all $i \in [1, m]$, then it would follow that
\begin{equation}\label{stationary2}
\m{x}^* = \arg \min \{ \C{L}(\m{x}, \g{\lambda}^*) : \m{x} \in \mathbb{R}^n \}.
\end{equation}
Since $\m{Ax}^* = \m{b}$, we conclude that
$\m{x}^*$ is an optimal solution of (\ref{Prob})--(\ref{ProbM}), and
$\g{\lambda}^*$ is an associated multiplier.
To show that (\ref{stationary1}) holds for all $i$, we need to consider
the situation where $\Gamma_i^k$ tends to infinity.

{\bf Case 2.} Suppose that $\Gamma_i^k$ approaches infinity.
Let $\bar{\m{x}}_i^k$ be the minimizer of $L_i^k$ defined in (\ref{lik}).
Observe that minimizing
$L_i^k (\m{u})$ over $\m{u} \in \mathbb{R}^{n_i}$ is equivalent to minimizing
a sum of the form $g(\m{u}) + h(\m{u}) + \langle \m{u}, \m{c}^k \rangle$ where
$h$ corresponds to $h_i$, $\m{c}^k = \m{A}_i\tr(\g{\lambda}^k - \rho\m{b}_i^k)$,
and $g(\m{u}) = f_i(\m{u}) + 0.5\rho \|\m{A}_i\m{u}\|^2$.
Note that $g$ is smooth and satisfies a strong convexity condition
\begin{equation}\label{strong}
(\m{u} - \m{v})\tr (\nabla g(\m{u}) - \nabla g(\m{v})) \ge
\rho \nu_i \|\m{u} - \m{v}\|^2,
\end{equation}
where $\nu_i > 0$ is the smallest eigenvalue of $\m{A}_i\tr \m{A}_i$.
By the strong convexity of $L_i^k$, it has a unique minimizer, and
from the first-order optimality conditions and the strong convexity
condition (\ref{strong}), we obtain the bound
\begin{equation}\label{lipbound}
\|\bar{\m{x}}_i^j - \bar{\m{x}}_i^k\| \le \|\m{c}^j - \m{c}^k\|/(\rho \nu_i).
\end{equation}
Since $\m{z}^k$, $\m{y}_+^k$, and $\g{\lambda}^k$ are bounded sequences,
it follows that $\bar{\m{x}}_i^k$ is a bounded sequence.
For $k \in \C{K}$, the sequences
$\m{z}^k$, $\m{y}_+^k$, and $\g{\lambda}^k$ converge to
$\m{x}^*$, $\m{x}_+^*$, and $\g{\lambda}^*$ respectively, which implies that
\begin{equation}\label{*}
\m{c}^* = \lim_{k\in\C{K}} \m{c}^k =
\m{A}_i\tr \left[ \g{\lambda}^* - \rho \left( \m{b}
- \sum_{j\ne i} \m{A}_j \m{x}_j^* \right) \right] =
\m{A}_i\tr \left[ \g{\lambda}^* - \rho \m{A}_i\m{x}_i^* \right] ,
\end{equation}
where the last equality is due to the identity $\m{Ax}^* = \m{b}$.
Consequently, by (\ref{lipbound}),
$\bar{\m{x}}_i^k$ for $k \in \C{K}$ forms a Cauchy
sequence which approaches a limit.
We use the stopping condition to determine the limit.

Let us insert $l = l_i^k$ and $\m{u}_i^l = \m{x}_i^{k+1}$ in the
inequality (\ref{xbarbound}).
By the linear independence of the columns of $\m{A}_i$ and the upper bound
(\ref{delta_bound}) for $\delta^l$, there exists $\beta > 0$ such that
\begin{eqnarray*}
\beta \|\bar{\m{x}}_i^k - \m{x}_i^{k+1}\|^2 &\le&
\frac{\rho}{\delta^l} \| \m{A}_i (\bar{\m{x}}_i^k - \m{u}_{i,k}^l)\|^2 \le
\|\bar{\m{x}}_i^k - \m{u}_{i,k}^{l-1}\|^2 -
\|\bar{\m{x}}_i^k - \m{u}_{i,k}^l\|^2 \\
&=& 2 \langle \bar{\m{x}}_i^k - \m{x}_i^{k+1} ,
\m{u}_{i,k}^{l} - \m{u}_{i,k}^{l-1} \rangle
+ \| \m{u}_i^{l} - \m{u}_i^{l-1}\|^2 \\
&\le& 2 \|\bar{\m{x}}_i^k - \m{x}_i^{k+1}\|
\|\m{u}_{i,k}^{l} - \m{u}_{i,k}^{l-1}\| +
\|\m{u}_{i,k}^{l} - \m{u}_{i,k}^{l-1}\|^2.
\end{eqnarray*}
We move the $\|\bar{\m{x}}_i^k - \m{x}_i^{k+1}\|$ term on the right
side of this inequality to the left side and complete the square to obtain
the relation
\[
\|\bar{\m{x}}_i^k - \m{x}_i^{k+1}\| \le
\frac{\|\m{u}_{i,k}^{l} - \m{u}_{i,k}^{l-1}\|}{\beta}
\left( 1 + \sqrt{\beta + 1} \right).
\]
Square this inequality and divide by $\Gamma_i^k$ to get
\[
\frac{\|\bar{\m{x}}_i^k - \m{x}_i^{k+1}\|^2}{\Gamma_i^k} \le
\frac{\|\m{u}_{i,k}^{l} - \m{u}_{i,k}^{l-1}\|^2}{\Gamma_i^k}
\left( \frac{1 + \sqrt{\beta + 1}}{\beta^2} \right)^2 .
\]
Since $l = l_i^k$, it follows from the stopping condition of Step~1c and
from (\ref{eklim}) that the right of this inequality approaches zero as
$k$ tends to infinity.
Earlier we showed that $\bar{\m{x}}_i^k$ is a bounded sequence.
Since $\Gamma_i^k$ tends to infinity in Case~2, and
$\bar{\m{x}}_i^k/\sqrt{\Gamma_i^k}$ approaches zero, we conclude that
$\m{x}_i^{k+1}/\sqrt{\Gamma_i^k}$ approaches zero.
Due to the inequality $\Gamma_i^{k+1} \ge \Gamma_i^k$,
$\m{x}_i^{k+1}/\sqrt{\Gamma_i^{k+1}}$ also approaches zero as $k$
tends to infinity.
Since $\|\m{x}_i^k -\bar{\m{x}}_i^k\| \le$
$\|\m{x}_i^k\| + \|\bar{\m{x}}_i^k\|$,
the right side of (\ref{GD-Converge}) approaches zero.
Hence, (\ref{GD-Converge}) implies that $\m{z}_i^k$ approaches
$\bar{\m{x}}_i^k$ as $k$ tends to infinity.
And since $\m{z}_i^k$ also approaches $\m{x}_i^*$ for $k \in \C{K}$,
we conclude that $\bar{\m{x}}_i^k$ approaches $\m{x}_i^*$ as
$k \in \C{K}$ tends to infinity.
Let $\bar{\m{x}}_i^*$ be defined by
\[
\bar{\m{x}}_i^* =
\arg \min_{\m{u}} \{ g(\m{u}) + h(\m{u}) +
\langle \m{u}, \m{c}^* \rangle \}.
\]
By (\ref{lipbound}) and the fact that
$\bar{\m{x}}_i^k$ approaches $\m{x}_i^*$ as
$k \in \C{K}$ tends to infinity, we conclude that
$\bar{\m{x}}_i^* = \m{x}_i^*$.
In summary, we have
\begin{eqnarray}
\lim_{k \in \C{K}} \bar{\m{x}}_i^k &=& \m{x}_i^* = \bar{\m{x}}_i^* =
\arg \min_{\m{u}} \{ g(\m{u}) + h(\m{u}) +
\langle \m{u}, \m{c}^* \rangle \}. \nonumber \\
&=& \arg \min_{\m{u}} \{ f_i(\m{u}) + 0.5\rho \|\m{A}_i \m{u}\|^2 + h_i (\m{u})
+ \langle \g{\lambda}^* - \rho \m{A}_i \m{x}_i^*, \m{u} \rangle \} .
\label{infinity}
\end{eqnarray}
The first-order optimality conditions for (\ref{infinity}) are
exactly the same as the first-order optimality conditions for
(\ref{stationary1}).
This shows that (\ref{stationary1}) holds in either Case~1 or Case~2.
Hence, (\ref{stationary2}) holds and 
$\m{x}^*$ is an optimal solution of (\ref{Prob})--(\ref{ProbM})
with associated multiplier $\g{\lambda}^*$.

Finally, we need to show that the entire sequence converges.
If $\Gamma_i^k$ is uniformly bounded as in Case~1, then by (\ref{uiklim}),
$\m{x}_i^k$ approaches $\m{x}_i^*$ and
$\|\m{x}_i^k - \m{x}_i^*\|^2/\Gamma_i^k$ approaches zero
as $k$ tends to infinity with $k \in \C{K}$.
On the other hand, when
$\Gamma_i^k$ tends to infinity as in Case~2, we showed that
$\|\m{x}_i^k - \bar{\m{x}}_i^k\|^2/\Gamma_i^k$ approaches zero.
Since $\bar{\m{x}}_i^k$ for $k \in \C{K}$ approaches $\m{x}_i^*$ by
(\ref{infinity}) and $\Gamma_i^k$ tends to infinity, it follows that
$\|\m{x}_i^k - \m{x}_i^*\|^2/\Gamma_i^k$ approaches zero too.
Thus in either Case~1 or Case~2,
$\|\m{x}_i^k - \bar{\m{x}}_i^k\|^2/\Gamma_i^k$ approaches zero as $k$ tends
to infinity with $k \in \C{K}$.
Letting $j$ tend to infinity in (\ref{Ej2}) with $j \in \C{K}$,
it follows that $E_j$ approaches zero.
Moreover, (\ref{Ej2}) implies that along the entire sequence,
$\m{y}_+^k$ approaches $\m{x}_+^*$
and $\g{\lambda}^k$ approaches $\g{\lambda}^*$.
By (\ref{L-lim12}), the entire sequence of iterates
$\m{z}_+^k$ approaches $\m{x}_+^*$.
Since $\m{Az}^k$ approaches $\m{b}$ (see (\ref{L-lim12})),
$\m{Ax}^* = \m{b}$, and $\m{A}_1\tr\m{A}_1$
is invertible, the entire sequence $\m{z}_1^k$ approaches $\m{x}_1^*$.
Finally, since $\m{y}_1^{k+1} = \m{z}_1^k$, we deduce that the entire
$\m{y}^k$ sequence approaches $\m{x}^*$.
This completes the proof.
\end{proof}
\smallskip

\section{Accelerated BOSVS}
\label{aBOSVS}
The inner loop for the accelerated BOSVS algorithm appears
in Algorithm~\ref{3}.
As we will see, the inner loop (Step~1) of accelerated BOSVS
converges to the minimizer of $L_i^k$, exactly as in multistep BOSVS;
however, the convergence speed of the multistep BOSVS inner loop is
$O(1/\sqrt{l})$ for the $\m{z}_i^k$ iterates and $O(1/l)$ for the objective
(see Remark~\ref{remark1}),
while the convergence speed in accelerated BOSVS is $O(1/l)$ for
the $\m{z}_i^k$ iterates and $O(1/l^2)$ for the objective,
which is optimal for first-order methods applied to general convex,
possibly nonsmooth optimization problems.
\renewcommand\figurename{Alg.}
\begin{figure}[h]
{\tt
\begin{tabular}{l}
\hline
{\bf Initialize:} $\m{a}_i^0 = \m{u}_i^0 = \m{x}_i^k$ and $\alpha^1 = 1$.
\\[.05in]
{\bf For } $l = 1, 2, \ldots $\\[.05in]
\begin{tabular}{ll}
1a. & Choose $\delta^l \ge \delta_{\min}$ and when $l > 1$, choose
$\alpha^l \in (0, 1)$ such that\\[.05in]
& $\quad f_i(\bar{\m{a}}_i^{l}) +$
$\langle \nabla f_i(\bar{\m{a}}_i^{l}), \m{a}_i^{l} -\bar{\m{a}}_i^{l-1} \rangle
+ \frac{(1-\sigma)\delta^l}{2\alpha^l} \|\m{a}_i^{l} - \bar{\m{a}}_i^{l-1}\|^2
\ge f_i(\m{a}_i^{l})$,\\[.05in]
& where $\m{a}_i^l = (1-\alpha^l)\m{a}_i^{l-1} + \alpha^l \m{u}_i^l$,
$\bar{\m{a}}_i^l =
(1-\alpha^l)\m{a}_i^{l-1} + \alpha^l \m{u}_i^{l-1}$, and \\[.05in]
&$\m{u}_i^l =
\arg \min \{ Q(\m{u}) + h_i (\m{u}): \m{u} \in \mathbb{R}^{n_i} \}$
with \\[.05in]
& \quad $Q(\m{u}) = 
\langle \nabla f_i(\bar{\m{a}}_i^{l}), \m{u} \rangle
+ \frac{\delta^l}{2} \|\m{u} - \m{u}_i^{l-1}\|^2 +
\frac{\rho}{2}\|\m{A}_i\m{u} - \bar{\m{a}}_i^k + \g{\lambda}^k/\rho\|^2$.
\\[.05in]
1b. &
If $\gamma^l =$
$(1/\delta^1) \displaystyle{\prod_{j=2}^l} (1-\alpha^j)^{-1}$
$\ge \Gamma_i^{k-1} $, where $\gamma^1 = 1/\delta^1$,\\[.05in]
& and $\| \m{a}_i^l - \m{a}_i^{l-1} \| \le \psi (e^{k-1})$, then break.
\end{tabular}\\
{\bf Next} \\[.05in]
1c. Set $\m{x}_i^{k+1}= \m{u}_i^l$,  $\m{z}_i^k = \m{a}_i^l$,
$\Gamma_i^k = \gamma^l$, and
$r_i^k = (1/\Gamma_i^k) \sum_{j=1}^l \|\m{u}_i^j - \m{u}_i^{j-1}\|^2$.\\[.05in]
\hline
\end{tabular}
}
\caption{Inner loop in Step~$1$ of Algorithm~$\ref{ADMMcommon}$
for the accelerated BOSVS scheme.}
\label{3} 
\end{figure}
\renewcommand\figurename{Fig.}

Two parameter sequences appear in the accelerated BOSVS scheme,
the $\delta^l$ and $\alpha^l$ sequences.
They must be chosen so that the line search condition of Step~1a is satisfied
for each value of $l$, and the stopping condition of Step~1b is satisfied
for $l$ sufficiently large.
If the Lipschitz constant $\zeta_i$ of $f_i$ is known, then we could take
\begin{equation}\label{AG_constant}
\delta^l = \frac{1}{(1-\sigma)}\frac{2 \zeta_i}{l } \quad \mbox{and}
\quad \alpha^l = \frac{2}{l+1} \in (0,1],
\end{equation}
in which case, we have
\[
\frac{(1-\sigma) \delta^l}{ \alpha^l} = \frac{(l+1) \zeta_i}{l} > \zeta_i.
\]
This relation along with a Taylor series expansion of $f_i$
around $\m{u}_i^{l-1}$ implies that the
line search condition in Step~1a of accelerated BOSVS is satisfied for each $l$.
Moreover, we show (after Lemma~\ref{lem-AG-Conv}) that with these choices
for $\delta^l$ and $\alpha^l$, the stopping condition of Step~1b is
also satisfied eventually.

A different, adaptive way to choose the parameters,
that does not require knowledge
of the Lipschitz constant for $f_i$, is the following:
Choose $\delta_0^l \in [\delta_{\min}, \delta_{\max}]$, where
$0 < \delta_{\min} < \delta_{\max} < \infty$ are safeguard parameters, and set
\begin{eqnarray}
\delta^l &=&
\frac{2}{\theta^l + \sqrt{(\theta^l)^2 +
4 \theta^l \Lambda^{l-1}}} \quad \mbox{and}
\quad \alpha^l = \frac{1}{1 + \delta^l \Lambda^{l-1}}, \quad \mbox{where}
\label{AG_linesearch}\\
\Lambda^l &=& \sum_{i=1}^l 1/\delta^i, \quad \Lambda^0 = 0,
\quad \mbox{and}\quad
\theta^l = 1/(\delta_0^l \eta^j) \mbox{ with } \eta > 1.
\nonumber
\end{eqnarray}
After some algebra, it can be shown that
\begin{equation}\label{delta/alpha}
\frac{\delta^l}{\alpha^l} = \frac{1}{\theta^l} = \delta_0^l \eta^j.
\end{equation}
Hence, the ratio $\delta^l/\alpha^l$ appearing in the line search condition
of Step~1a tends to infinity as $j$ tends to infinity since $\eta > 1$.
We take $j \ge 0$ to be the smallest integer for which the line search
condition is satisfied.
Based on the identity (\ref{delta/alpha}), the expression
$\delta^l/\alpha^l$ has exactly the same effect as $\delta_i^k$ in
generalized BOSVS.
Consequently, it satisfies exactly the same inequality (\ref{delta_bound}).

Let us first observe that when
$e^k =0$, we have reached a solution of (\ref{Prob})--(\ref{ProbM}).
\smallskip

\begin{lemma}\label{L-stop-cond3}
If $e^{k}=0$ in the accelerated BOSVS algorithm,
then $\m{x}^{k+1} = \m{x}^k = \m{y}^k = \m{z}^k$ solves
$(\ref{Prob})$--$(\ref{ProbM})$ and $(\m{x}^k, \g{\lambda}^k) \in \C{W}^*$.
\end{lemma}
\smallskip

\begin{proof}
If $e^k = 0$, then $r_i = 0$ for each $i$.
It follows that
\begin{equation}\label{equals}
\m{x}_i^k = \m{u}_i^0 = \m{u}_i^1 = \ldots = \m{u}_i^l.
\end{equation}
By Step~1c, $\m{u}_i^l = \m{x}_i^{k+1}$.
By the definitions
$\m{a}_i^l = (1-\alpha^l)\m{a}_i^{l-1} + \alpha^l \m{u}_i^l$ and
$\bar{\m{a}}_i^l = (1-\alpha^l)\m{a}_i^{l-1} + \alpha^l \m{u}_i^{l-1}$
where $\m{u}_i^0 = \m{x}_i^k$, we have
$\m{a}_i^l = \bar{\m{a}}_i^l = \m{x}_i^k$ for each $l$ due to (\ref{equals}).
Again, by Step~1c, $\m{z}_i^k = \m{x}_i^k$.
Consequently, we have $\m{x}^{k+1} = \m{x}^k = \m{z}^k$.
Since all three algorithms in this paper share
Algorithm~\ref{ADMMcommon}, the remainder of the proof is
exactly as in Lemma~\ref{L-stop-cond}.
\end{proof}

We now establish the following analogue of Lemma~\ref{lem-GD-Conv}.
\begin{lemma}\label{lem-AG-Conv}
If the inner loop sequence $\xi^l := \delta^l \alpha^l \gamma^l$
associated with accelerated BOSVS is nonincreasing as a function of $l$,
then for each $i \in [1,m]$, we have
\begin{equation}\label{AG-Converge}
\quad \quad
\nu_i \rho \| \m{z}_i^{k} - \bar{\m{x}}_i^k\|^2 
+ \frac{\sigma}{\Gamma_i^{k}}
\sum_{l=1}^{l_i^k}  \xi^l \|\m{u}_i^l- \m{u}_i^{l-1}\|^2 
\le  \frac{\| \m{x}_i^k - \bar{\m{x}}_i^k\|^2}
{\Gamma_i^{k}} ,
\end{equation}
where $l_i^k$ is the terminating value of $l$ at iteration $k$,
$\bar{\m{x}}_i^k$ is the minimizer of the function $L_i^k$ defined
in $(\ref{lik})$, and $\nu_i>0$ is the smallest eigenvalue of
$\m{A}_i \tr \m{A}_i$.
\end{lemma}
\smallskip

\begin{proof}
By the definition
$\m{a}_i^l = (1-\alpha^l) \m{a}_i^{l-1} + \alpha^l \m{u}_i^l$, we have
\[
\langle \nabla  f_i({\bar{\m{a}}}_i^l ),
\m{a}_i^l  - {\bar{\m{a}}}_i^l \rangle =
(1 - \alpha^l)
\langle \nabla  f_i({\bar{\m{a}}}_i^l) ,
\m{a}_i^{l-1}  - {\bar{\m{a}}}_i^l \rangle
+ \alpha^l \langle \nabla  f_i({\bar{\m{a}}}_i^l ),
\m{u}_i^l - {\bar{\m{a}}}_i^l \rangle .
\]
Add to this the identity
$f_i({\bar{\m{a}}}_i^l ) = (1-\alpha^l) f_i({\bar{\m{a}}}_i^l ) +
\alpha^l f_i({\bar{\m{a}}}_i^l )$ to obtain
\begin{eqnarray*}
&f_i({\bar{\m{a}}}_i^l ) +
\langle \nabla  f_i({\bar{\m{a}}}_i^l ),
\m{a}_i^l  - {\bar{\m{a}}}_i^l \rangle =& \nonumber \\[.05in]
&(1 - \alpha^l) \left[ f_i({\bar{\m{a}}}_i^l ) +
\langle \nabla  f_i({\bar{\m{a}}}_i^l) ,
\m{a}_i^{l-1}  - {\bar{\m{a}}}_i^l \rangle \right]
+ \alpha^l \left[ f_i({\bar{\m{a}}}_i^l )
+ \langle \nabla  f_i({\bar{\m{a}}}_i^l ),
\m{u}_i^l - {\bar{\m{a}}}_i^l \rangle \right].&
\end{eqnarray*}
By the convexity of $f_i$, it follows that
\[
f_i({\bar{\m{a}}}_i^l ) +
\langle \nabla  f_i({\bar{\m{a}}}_i^l) ,
\m{a}_i^{l-1}  - {\bar{\m{a}}}_i^l \rangle \le f_i(\m{a}_i^{l-1}) .
\]
Hence, we have
\[
f_i({\bar{\m{a}}}_i^l ) +
\langle \nabla  f_i({\bar{\m{a}}}_i^l ) ,
\m{a}_i^l  - {\bar{\m{a}}}_i^l \rangle  \le
(1 - \alpha^l) f_i(\m{a}_i^{l-1}) + 
\alpha^l \left[ f_i({\bar{\m{a}}}_i^l )
+ \langle \nabla  f_i({\bar{\m{a}}}_i^l ),
\m{u}_i^l - {\bar{\m{a}}}_i^l \rangle \right] .
\]
Adding and subtracting any $\m{u} \in \mathbb{R}^{n_i}$ in the last
term, and then exploiting the convexity of $f_i$ gives
\begin{eqnarray*}
f_i({\bar{\m{a}}}_i^l )
+ \langle \nabla  f_i({\bar{\m{a}}}_i^l ),
\m{u}_i^l - {\bar{\m{a}}}_i^l \rangle &=&
\left[ f_i({\bar{\m{a}}}_i^l )
+ \langle \nabla  f_i({\bar{\m{a}}}_i^l ),
\m{u} - {\bar{\m{a}}}_i^l \rangle \right] +
\langle \nabla  f_i({\bar{\m{a}}}_i^l ),
\m{u}_i^l - \m{u} \rangle \\
&\le& f_i(\m{u}) + \langle \nabla  f_i({\bar{\m{a}}}_i^l ),
\m{u}_i^l - \m{u} \rangle .
\end{eqnarray*}
Therefore,
\begin{equation} \label{1234}
\quad \quad f_i({\bar{\m{a}}}_i^l ) +
\langle \nabla  f_i({\bar{\m{a}}}_i^l ),
\m{a}_i^l  - {\bar{\m{a}}}_i^l \rangle \le
(1 - \alpha^l) f_i (\m{a}_i^{l-1} ) +
\alpha^l [ f_i(\m{u}) + \langle \nabla
f_i({\bar{\m{a}}}_i^l ), \m{u}_i^l  - \m{u} \rangle ].
\end{equation}

Now by the line search condition in Step~1a of accelerated BOSVS and
then by (\ref{1234}), we have
\begin{eqnarray*}
L_i^k (\m{a}_i^l) &=& f_i (\m{a}_i^l) +
\frac{\rho}{2} \|\m{A}_i \m{a}_i^l - \m{b}_i^k+
\g{\lambda}^k/\rho\|^2 + h_i  (\m{a}_i^l) \nonumber \\
&\le& f_i({\bar{\m{a}}}_i^l) +
\langle \nabla f_i({\bar{\m{a}}}_i^l), \m{a}_i^l -{\bar{\m{a}}}_i^l \rangle +
\frac{(1-\sigma)\delta^l}{2 \alpha^l} \|\m{a}_i^l - {\bar{\m{a}}}_i^l\|^2
\nonumber
\\[.05in]
&&\quad \quad \quad +
\frac{\rho}{2} \|\m{A}_i \m{a}_i^l - \m{b}_i^k + \g{\lambda}^k/\rho\|^2 +
h_i  (\m{a}_i^l) \nonumber \\[.05in]
&\le&
(1 - \alpha^l) f_i (\m{a}_i^{l-1} ) +
\alpha^l f_i(\m{u}) + \alpha^l \langle \nabla
f_i({\bar{\m{a}}}_i^l ), \m{u}_i^l  - \m{u} \rangle +
\frac{(1-\sigma)\delta^l}{2 \alpha^l} \|\m{a}_i^l - {\bar{\m{a}}}_i^l\|^2
\nonumber
\\[.05in]
&& \quad \quad \quad +
\frac{\rho}{2} \|\m{A}_i \m{a}_i^l - \m{b}_i^k + \g{\lambda}^k/\rho\|^2 +
h_i  (\m{a}_i^l). \label{1235}
\end{eqnarray*}
Next, we utilize the definitions of $\m{a}_i^l$ and $\bar{\m{a}}_i^l$ and
the convexity of both $h_i$ and the norm term to obtain
\begin{eqnarray}
L_i^k (\m{a}_i^l) &\le&
(1 - \alpha^l) f_i (\m{a}_i^{l-1} ) +
\alpha^l [ f_i(\m{u}) +  \langle \nabla  f_i({\bar{\m{a}}}_i^l ), 
\m{u}_i^l  - \m{u} \rangle] +
\frac{(1-\sigma)\delta^l}{2 \alpha^l} \|\m{a}_i^l - {\bar{\m{a}}}_i^l\|^2
\nonumber \\
&&  \quad + (1-\alpha^l) \left(\frac{\rho}{2}
\|\m{A}_i \m{a}_i^{l-1} - \m{b}_i^k  + \g{\lambda}^k/\rho \|^2 
+ h_i  (\m{a}_i^{l-1}) \right) \nonumber \\
&& \quad + \alpha^l \left(\frac{ \rho}{2} \|\m{A}_i \m{u}_i^l - \m{b}_i^k
+ \g{\lambda}^k/\rho \|^2 + h_i  (\m{u}_i^l) \right) \nonumber\\
&=&  (1-\alpha^l) \left( f_i (\m{a}_i^{l-1} ) +
\frac{\rho}{2} \|\m{A}_i \m{a}_i^{l-1} - \m{b}_i^k +
\g{\lambda}^k/\rho \|^2 +  h_i  (\m{a}_i^{l-1})  \right) \nonumber\\ 
&&\quad + \alpha^l [ f_i(\m{u}) + \langle \nabla
 f_i({\bar{\m{a}}}_i^l ),\m{u}_i^l  - \m{u} \rangle]
+ \frac{(1-\sigma)\delta^l \alpha^l}{2}
\|\m{u}_i^l - \m{u}_i^{l-1}\|^2   \nonumber\\ 
&&\quad + \alpha^l \left(\frac{ \rho}{2}
\|\m{A}_i \m{u}_i^l - \m{b}_i^k +
\g{\lambda}^k/\rho \|^2 + h_i  (\m{u}_i^l) \right) \nonumber \\
&=&  (1-\alpha^l) L_i^k (\m{a}_i^{l-1} ) 
+ \alpha^l [ f_i(\m{u}) +
\langle \nabla  f_i({\bar{\m{a}}}_i^l ),\m{u}_i^l  - \m{u} \rangle]
\nonumber \\
&& \quad
+ \frac{(1-\sigma)\delta^l \alpha^l}{2 } \|\m{u}_i^l - \m{u}_i^{l-1}\|^2
+  \alpha^l \left(\frac{ \rho}{2} \|\m{A}_i \m{u}_i^l - \m{b}_i^k
+ \g{\lambda}^k/\rho\|^2 + h_i  (\m{u}_i^l) \right).  \label{qwer}
\end{eqnarray}
Since $h_i$ is convex, we have
\begin{equation}\label{hineq}
h_i(\m{u}^l) + \langle \m{p}, \m{u} - \m{u}^l \rangle \le h_i(\m{u})
\end{equation}
for any $\m{p} \in \partial h_i (\m{u}^l)$.
The expansion of the quadratic $Q$ in Step~1a of accelerated BOSVS
around $\m{u}^l$ can be written
\begin{equation}\label{Qeq}
Q(\m{u}^l) + \nabla Q(\m{u}^l)(\m{u} - \m{u}^l)
+ \frac{1}{2}(\m{u} - \m{u}^l)\tr (\delta^l \m{I} + \rho \m{A}_i\tr\m{A}_i)
(\m{u} - \m{u}^l) = Q(\m{u}).
\end{equation}
Since $\m{u}^l$ minimizes $Q + h_i$ in Step~1a, the first-order optimality
conditions imply that $\m{p} + \nabla Q(\m{u}^l) = \m{0}$ for some
$\m{p} \in \partial h_i (\m{u}^l)$.
We choose $\m{p} = -\nabla Q(\m{u}^l)$, and then
multiply (\ref{hineq}) and (\ref{Qeq}) by $\alpha^l$ and add to (\ref{qwer})
to obtain
\begin{eqnarray*}
L_i^k(\m{a}_i^l) &\le&
(1-\alpha^l) L_i^k (\m{a}_i^{l-1} ) + \alpha^l
L_i^k(\m{u})  + \frac{ \delta^l \alpha^l }{2}
( \|\m{u} - \m{u}_i^{l-1}\|^2 - \|\m{u} - \m{u}_i^l \|^2) \\
&&\quad -\frac{\sigma\delta^l \alpha^l}{2 } \|\m{u}_i^l - \m{u}_i^{l-1}\|^2  
 - \frac{\alpha^l \rho}{2}\|\m{A}_i  (\m{u} -\m{u}_i^l)\|^2.
\end{eqnarray*}
Hence, for any $\m{u} \in \mathbb{R}^{n_i}$ we have
\begin{eqnarray}
L_i^k (\m{a}_i^l) - L_i^k(\m{u})
&\le&(1-\alpha^l) ( L_i^k (\m{a}_i^{l-1} ) - L_i^k(\m{u}) )
+ \frac{ \delta^l \alpha^l }{2} ( \|\m{u} - \m{u}_i^{l-1}\|^2
- \|\m{u} - \m{u}_i^l \|^2)  \nonumber \\
&& \quad -\frac{\sigma\delta^l \alpha^l}{2 } \|\m{u}_i^l - \m{u}_i^{l-1}\|^2 
 - \frac{\alpha^l \rho}{2}\|\m{A}_i  (\m{u} -\m{u}_i^l)\|^2.
\label{asdf}
\end{eqnarray}

From the definition of $\gamma^l$ in accelerated BOSVS, it follows that
$(1-\alpha^l) \gamma^l = \gamma^{l-1}$ with the convention that
$\gamma^0 = 0$ (since $\alpha^1 = 1)$.
Hence, for any sequence $d^l$, $l \ge 0$, we have
\begin{equation}\label{id1}
\sum_{l=1}^j \left( \gamma^l d^l - (1-\alpha^l)\gamma^l d^{l-1}
\right) =
\sum_{l=1}^j \left( \gamma^l d^l - \gamma^{l-1} d^{l-1} \right) =
\gamma^j d^j.
\end{equation}
Suppose that $d^l \ge 0$ for each $l$.
By assumption, $\xi^l = \gamma^l \delta^l \alpha^l$ is nonincreasing;
since $\alpha^1 = 1$ and $\gamma^1 = 1/\delta^1$, it follows that
$\xi^1 = 1$, and we have
\begin{eqnarray}
\sum_{l=1}^j \xi^l \left( d^l - d^{l-1} \right) &=&
d^1 - d^0 + \sum_{l=2}^j \xi^l \left( d^l - d^{l-1} \right)
\label{id2} \\
&\ge& d^1 - d^0 + \sum_{l=2}^j \left( \xi^l d^l - \xi^{l-1} d^{l-1} \right) =
\xi^j d^j - d^0 . \nonumber
\end{eqnarray}
We now multiply (\ref{asdf}) by $\gamma^l$ and sum over
$l$ between 1 and $l_i^k$.
Exploiting the identity (\ref{id1}) with
$d^l = L_i^k (\m{a}_i^l) - L_i^k(\m{u})$ and
(\ref{id2}) with $d^l = \|\m{u}_i^l - \m{u}\|^2$, we obtain
\begin{eqnarray}
L_i^k(\m{u}) -  L_i^k (\m{a}_i^{l_i^k})
&\ge& \frac{1}{2 \Gamma_i^k}
(\xi^{l_i^k} \|\m{u} - \m{u}_i^{l_i^k} \|^2 -   \|\m{u} - \m{u}_i^{0}\|^2 )
+ \frac{\sigma}{2 \Gamma_i^k} \sum_{l=1}^{l_i^k} 
\xi^l  \|\m{u}_i^l - \m{u}_i^{l-1}\|^2  \nonumber\\
&& \quad +
\frac{\rho}{2 \Gamma_i^k}
\sum_{l=1}^{l_i^k} (\gamma^l \alpha^l) \|\m{A}_i  (\m{u} -\m{u}_i^l)\|^2,
\label{zxcv}
\end{eqnarray} 
where $\Gamma_i^k$ denotes the final $\gamma^l$ in accelerated BOSVS.

Next, we multiply the definition
$\m{a}_i^j = (1-\alpha^j) \m{a}_i^{j-1} + \alpha^l \m{u}_i^j$
by $\gamma^j$ and sum over $j$ between 1 and $l$.
Again, exploiting the identity
$(1-\alpha^j) \gamma^j = \gamma^{j-1}$ yields
\begin{equation}\label{convex1}
\m{a}_i^l =
\frac{1}{\gamma^l} \sum_{j=1}^{l} (\gamma^j \alpha^j) \m{u}_i^j. 
\end{equation}
Since $\alpha^j \gamma^j = \gamma^j - \gamma^{j-1}$, it follows that
\begin{equation}\label{convex2}
\gamma^l = \sum_{j=1}^l \alpha^j \gamma^j .
\end{equation}
Consequently, $\m{a}_i^l$ is a convex combination of
$\m{u}_i^1$ through $\m{u}_i^l$.
Since $\|\m{A}_i  (\m{u} - \m{w})\|^2$ is a convex function of $\m{w}$,
Jensen's inequality yields
\[
\frac{1}{\Gamma_i^k} \sum_{l=1}^{l_i^k} (\gamma^l \alpha^l)
\|\m{A}_i  (\m{u} - \m{u}_i^l)\|^2 
\ge \|\m{A}_i  (\m{u} - \m{a}_i^{l_i^k})\|^2 =
\|\m{A}_i  (\m{u} - \m{z}_i^k)\|^2 .
\]
We apply this inequality to the last term in (\ref{zxcv}) and substitute
$\m{z}_i^k = \m{a}_i^{l_i^k}$,
$\m{x}_i^{k+1} = \m{u}_i^{l_i^k}$, and $\m{x}_i^k = \m{u}_i^0$ to obtain
\begin{eqnarray}
L_i^k(\m{u}) -  L_i^k (\m{z}_i^{k})
&\ge& \frac{1}{2 \Gamma_i^k}
(\xi^{l_i^k} \|\m{u} - \m{x}_i^{k+1} \|^2 -   \|\m{u} - \m{x}_i^{k}\|^2 )
+ \frac{\sigma}{2 \Gamma_i^k} \sum_{l=1}^{l_i^k} 
\xi^l  \|\m{u}_i^l - \m{u}_i^{l-1}\|^2  \nonumber\\
&& \quad +
\frac{\rho}{2}
\|\m{A}_i  (\m{u} -\m{z}_i^k)\|^2.
\label{generalu}
\end{eqnarray} 
Finally, take $\m{u} = \bar{\m{x}}_i^k$.
Since the left side of (\ref{generalu}) is nonpositive for this choice
of $\m{u}$, the proof is complete.
\end{proof}

Let us now examine the assumptions and consequences of Lemma~\ref{lem-AG-Conv}
in the context of the choices (\ref{AG_constant}) and (\ref{AG_linesearch})
for the parameters $\delta^l$ and $\alpha^l$.
For the choice (\ref{AG_constant}) and for $l \ge 2$, we have
\begin{equation}\label{gammal1}
\gamma^l =
\frac{1}{\delta^1}\prod_{j=2}^l (1-\alpha^j)^{-1} =
\frac{1}{\delta^1}\prod_{j=2}^l \frac{j+1}{j-1} =
\frac{1}{\delta^1} \frac{l(l+1)}{2}.
\end{equation}
Hence, $\gamma^l$ is $O(l^2)$.
Since $\delta^l = \delta^1/l$, it follows that for $l \ge 2$,
\[
\xi^l := \delta^l \alpha^l \gamma^l =
\left( \frac{\delta^1}{l} \right)
\left( \frac{2}{l+1} \right) \left( \frac{l(l+1)}{2 \delta^1} \right) = 1.
\]
In the special case $l = 1$, $\xi^1 = \delta^1/\delta^1 = 1$.
Since the sequence $\xi^l$ is identically one, it is nonincreasing
and the assumption of Lemma~\ref{lem-AG-Conv} is satisfied.
Since $\Gamma_i^k$ is the final value for $\gamma^l$ in Step~1 of
accelerated BOSVS, it follows from (\ref{AG-Converge})
that $\|\m{z}_i^k - \bar{\m{x}}_i^k\| = O(1/l_i^k)$.

For the choice (\ref{AG_linesearch}) and for $l \ge 2$, we have
$\Lambda^l = (1/\delta^l) + \Lambda^{l-1}$ and
$\alpha^l = (1/\delta^l)/\Lambda^l$.
It follows that $1 - \alpha^l = \Lambda^{l-1}/\Lambda^l$ and for $l \ge 2$,
we have
\[
\gamma^l = \frac{1}{\delta^1} \prod_{j=2}^l (1-\alpha^j)^{-1}
=  \frac{1}{\delta^1} \prod_{j=2}^l (\Lambda^{j}/\Lambda^{j-1})
= \frac{1}{\delta^1} \frac{\Lambda^l}{\Lambda^1} = \Lambda^l .
\]
Hence,
\[
\xi^l := \delta^l \alpha^l \gamma^l =
\delta^l \left ( \frac{1/\delta^l}{\Lambda^l} \right) \Lambda^l  = 1.
\]
In the special case $l = 1$, we also have $\xi^1 = 1$.
Again, the sequence $\xi^l$ is identically one, which satisfies the
requirement of Lemma~\ref{lem-AG-Conv};
consequently, the speed with which $\m{z}_i^k$ converges to $\bar{\m{x}}_i^k$ 
depends on the growth rate of $\gamma^l$.
By the definition of $\gamma^l$ in accelerated BOSVS,
\begin{equation}\label{gammarecur}
\sqrt{\gamma^l} - \sqrt{\gamma^{l-1}} = 
\sqrt{\gamma^l} - \sqrt{(1- \alpha^l) \gamma^l} =
\left( 1- \sqrt{1-\alpha^l} \right)
\sqrt{\gamma^l} \ge \frac{\alpha^l \sqrt{\gamma^l}}{2}.
\end{equation}
Since $\xi^l := \delta^l \alpha^l \gamma^l = 1$,
it follows from (\ref{delta/alpha}) that
$(\alpha^l/\theta^l)\alpha^l\gamma^l =$
$(\alpha^l)^2\gamma^l/\theta^l = 1$,
which implies that
\begin{equation}\label{alphalgamma}
\alpha^l \sqrt{\gamma^l} = \sqrt{\theta^l}.
\end{equation}
By (\ref{gammarecur}), we have
\begin{equation}\label{recur}
\sqrt{\gamma^l} - \sqrt{\gamma^{l-1}} \ge \frac{\sqrt{\theta^l}}{2} .
\end{equation}
As noted beneath (\ref{delta/alpha}), $1/\theta^l$
satisfies the inequality (\ref{delta_bound}) for $\delta_i^k$,
which implies that
\begin{equation}\label{thetallower}
\theta^l \ge \Theta := \frac{1-\sigma}{\eta \zeta_i + (1-\sigma) \delta_{\max}}.
\end{equation}
Hence, (\ref{recur}) yields
\[
\sqrt{\gamma^l} - \sqrt{\gamma^{l-1}} \ge \frac{\sqrt{\Theta}}{2}.
\]
Since $\gamma^1 = 1/\delta^1 = \theta^1$, it follows that
\[
\sqrt{\gamma^l} \ge \sqrt{\Theta} +
\left( \frac{l-1}{2} \right) \sqrt{\Theta} \ge
\left( \frac{l}{2}\right)  \sqrt{\Theta} \quad \mbox{or} \quad
\gamma^l \ge \left( \frac{l^2}{4}\right)  \Theta,
\]
which implies that $\gamma^l = O(l^2)$.
In summary,
for either of the choices (\ref{AG_constant}) or (\ref{AG_linesearch}),
we have $\xi^l = 1$ for each $l$, and
$\|\m{z}_i^k - \bar{\m{x}}_i^k\| = O(1/l_i^k)$.
Moreover, by the inequality (\ref{generalu}) with
$\m{u} = \bar{\m{x}}_i^k$, the objective value satisfies
$L_i^k(\m{z}_i^k) - L_i^k(\bar{\m{x}}_i^k) = O(1/(l_i^k)^2)$.

Although Lemma~\ref{lem-AG-Conv} was stated in terms of the terminating
iteration $l_i^k$ of the inner iteration, it applies to any of the
inner iterations; that is, for each $i$ and $l$, we have
\[
\nu_i \rho \| \m{a}^{l} - \bar{\m{x}}_i^k\|^2 
+ \frac{\sigma}{\gamma^{l}}
\sum_{j=1}^{l}  \xi^j \|\m{u}_i^j- \m{u}_i^{j-1}\|^2 
\le  \frac{\| \m{x}_i^k - \bar{\m{x}}_i^k\|^2}
{\gamma^{l}} .
\]
Whenever $\gamma^l$ approaches infinity, as it does with the choices
(\ref{AG_constant}) and (\ref{AG_linesearch}), the right side approach
zero and $\m{a}^l$ converges to $\bar{\m{x}}_i^k$.
Hence, the stopping conditions in Step~1b of accelerated BOSVS are
satisfied for $l$ sufficiently large when $e^{k-1} \ne 0$.

The convergence of accelerated BOSVS, like the other algorithms,
relies on a decay property for the iterates, which we now give.
\begin{lemma}\label{L-key-lemma3}
If the accelerated BOSVS parameters
$\gamma^l$ tend infinity as $l$ grows and
$\xi^l := \delta^l \alpha^l \gamma^l = 1$ for each $l$,
then Lemma~$\ref{L-key-lemma2}$ holds for the accelerated scheme.
\end{lemma}
\smallskip

\begin{proof}
We substitute $\m{u} = \m{x}_i^*$ and $\xi^l = 1$ in (\ref{generalu}) to obtain
\[
L_i^k(\m{x}_i^*) - F_i^k(\m{z}_i^k) \ge
\frac{1}{2\Gamma_i^k} \left(
\|\m{x}_i^{k+1} - \m{x}_i^* \|^2 -   \|\m{x}_i^{k} - \m{x}_i^*\|^2 \right)
+ \frac{\sigma}{2 \Gamma_i^k} \sum_{l=1}^{l_i^k} 
\|\m{u}_i^l - \m{u}_i^{l-1}\|^2 ,
\]
where $F_i^k(\m{w}) = L_i^k(\m{w}) + (\rho/2) \|\m{A}_i(\m{w} - \m{x}_i^*)\|^2$.
This is exactly the same as (\ref{jumppoint}) in the proof of
Lemma~\ref{L-key-lemma2}.
The remainder of the proof is exactly as in the proof of
Lemma~\ref{L-key-lemma2}.
\end{proof}

Using the decay property of Lemmas~\ref{L-key-lemma2} and \ref{L-key-lemma3},
we now obtain the convergence of accelerated BOSVS.
\smallskip

\begin{theorem}\label{L-glob-thm3}
Suppose that for the inner loop sequence $\xi^l := \delta^l \alpha^l \gamma^l$
associated with accelerated BOSVS we have $\xi^l = 1$ for each $l$,
$\gamma^l$ tends to infinity as $l$ grows,
and there exists a constant $\kappa > 0$ such that
$\gamma^l (\alpha^l)^2 \ge \kappa$ for all $l$.
If accelerated BOSVS performs an infinite number of iterations
generating iterates $\m{y}^k$, $\m{z}^k$, and $\g{\lambda}^k$,
then the sequences
$\m{y}^k$ and $\m{z}^k$ both approach a common limit $\m{x}^*$ and
$\g{\lambda}^k$ approaches a limit $\g{\lambda}^*$ where
$(\m{x}^*, \g{\lambda}^*) \in \C{W}^*$.
\end{theorem}
\smallskip

\begin{proof}
The proof is identical to that of Theorem~\ref{L-glob-thm2} through the
end of Case~1.
For accelerated BOSVS, the fact that $\m{z}_i^k$ is a convex combination
of $\m{u}_{i,k}^l$ is shown in (\ref{convex1})--(\ref{convex2}).
The treatment of accelerated BOSVS first differs from that of multistep BOSVS
in the second paragraph of Case~2 ($\Gamma_i^k$ tends to $+\infty$)
where the multistep BOSVS stopping condition
$\| \m{u}_i^l - \m{u}_i^{l-1} \|/\sqrt{\gamma^l} \le \psi (e^{k-1})$,
is used to show that
\break
$\|\m{x}_i^k - \bar{\m{x}}_i^k\|^2/\Gamma_i^k$
approaches zero.
Since accelerated BOSVS uses the new stopping condition
$\|\m{a}_i^l - \m{a}_i^{l-1} \| \le \psi (e^{k-1})$,
a new analysis is needed in Case~2.

By the definition of $\m{a}^l$, we have
\[
\|\m{a}^l - \m{a}^{l-1}\| = \alpha^l \|\m{u}^l - \m{a}^{l-1}\| \ge
\alpha^l (\|\m{u}^l - \m{a}^{l}\| - \|\m{a}^l - \m{a}^{l-1}\|).
\]
If $\psi_k$ denotes $\psi(e^{k-1})$ and $l = l_i^k$ so that $\m{a}^l$ satisfies
the stopping criterion $\|\m{a}_i^l - \m{a}_i^{l-1}\| \le \psi_k$, then
\[
\alpha^l \|\m{u}^l - \m{a}^{l}\| \le
(1+\alpha^l) \|\m{a}^l - \m{a}^{l-1}\| \le 2\psi_k
\quad \mbox{or}
\quad \|\m{x}_i^{k+1} - \m{z}_i^{k}\| \le \frac{2\psi_k}{\alpha^l}
\]
since $\m{u}^l = \m{x}_i^{k+1}$ and $\m{a}^l = \m{z}_i^k$ when $l = l_i^k$.
Squaring this, dividing by $\gamma^l = \Gamma_i^k$,
and utilizing the assumption
that $\gamma^l (\alpha^l)^2 \ge \kappa$ for all $l$, we deduce that
\begin{equation}\label{psik}
\frac{\|\m{x}_i^{k+1} - \m{z}_i^k\|^2}{\Gamma_i^k} \le
\frac{4\psi_k^2}{\kappa}.
\end{equation}
Since $\psi_k$ approach zero by (\ref{eklim}), it follows that
$\|\m{x}_i^{k+1} - \m{z}_i^k\|^2/\Gamma_i^k$ approaches zero as $k$
tends to infinity.
Since $\Gamma_i^k$ is nondecreasing,
$\|\m{x}_i^{k+1} - \m{z}_i^k\|^2/\Gamma_i^{k+1}$ also approaches zero as $k$
tends to infinity.
Since $\m{z}_i^k$ is a bounded sequence and
$\Gamma_i^k$ tends to infinity in Case~2, we can replace
$\m{z}_i^k$ by any other bounded sequence and reach the same conclusion.
In particular, since the sequence $\bar{\m{x}}_i^k$ is bounded we conclude that
$\|\m{x}_i^{k} - \bar{\m{x}}_i^k\|^2/\Gamma_i^k$ approaches zero as $k$
tends to infinity, the same conclusion we reached in
multistep BOSVS scheme.
The rest of the proof is exactly as in Theorem~\ref{L-glob-thm2}.
This completes the proof.
\end{proof}

\begin{remark}
The parameter choices given in both $(\ref{AG_constant})$ and
$(\ref{AG_linesearch})$ satisfy the assumption of Theorem~$\ref{L-glob-thm3}$
that $\gamma^l (\alpha^l)^2 \ge \kappa > 0$ for some constant $\kappa$.
In particular, for $(\ref{AG_constant})$, we show in $(\ref{gammal1})$ that
$\gamma^l = l(l+1)/(2\delta^1)$.
This is combined with the definition of $\alpha^l$ in $(\ref{AG_constant})$ to
obtain
\[
\gamma^l (\alpha^l)^2 = \frac{2l}{\delta^1 (l+1)} \ge \frac{1}{\delta^1}
\]
for $l \ge 1$.
For the choice $(\ref{AG_linesearch})$, it follows from
$(\ref{alphalgamma})$ and $(\ref{thetallower})$ that
\[
\gamma^l (\alpha^l)^2 \ge
\Theta := \frac{1-\sigma}{\eta \zeta_i + (1-\sigma) \delta_{\max}}.
\]
\end{remark}

\begin{remark}
In this paper, we have focused on algorithms based on an inexact minimization
of $L_i^k$ in Step~$1$ of Algorithm~$\ref{ADMMcommon}$.
In cases where $f_i$ and $h_i$ are simple enough that the exact minimizer
$\bar{\m{x}}_i^k$ of $L_i^k$ can be quickly evaluated, we could simply set
$\m{x}_i^{k+1} = \m{z}_i^k = \bar{\m{x}}_i^k$ and $r_i^k = 0$ in Step~$1$.
The analysis of this exact algorithm is very similar to the analysis in
Theorems~$\ref{L-glob-thm2}$ and $\ref{L-glob-thm3}$.
In the analysis of the inexact algorithms, a key inequality
$(\ref{combine})$ was
\[
- \rho \left\langle \sum_{j \le i} \m{A}_j \m{z}_{e,j}^{k} +
\sum_{j > i} \m{A}_j \m{y}_{e,j}^k
+ \g{\lambda}_e^k/\rho, \; \m{A}_i \m{z}_{e,i}^{k} \right\rangle \ge \tau_i^k,
\]
where
\[
\tau_i^k = \frac{1}{2\Gamma_i^k} (\|\m{x}_{e,i}^{k+1}\|^2 - \|\m{x}_{e,i}^k\|^2)
+ \frac{\sigma}{2\Gamma_i^k}
\sum_{l = 1}^{l_i^k} \|\m{u}_{i,k}^{l} - \m{u}_{i,k}^{l-1}\|^2 .
\]
For an exact minimizer of $L_i^k$, the same inequality can be established
but with $\tau_i^k$ replaced by zero.
This follows directly from the first-order optimality conditions for a
minimizer of $L_i^k$ and for a minimizer of {\rm (\ref{Prob})--(\ref{ProbM})}.
Since $\tau_i^k$ disappears, then so do the $\m{x}_i^k$ and
$\m{u}_{i,k}^l$ terms in Lemma~{\rm \ref{L-key-lemma2}};
consequently, the analysis becomes simpler when the minimizer of
$L_i^k$ is exact.
\end{remark}
\section{Numerical Experiments}
\label{numerical}
In this section, we investigate the performance of the algorithms for
an image reconstructed problem that can be formulated as
\begin{equation}\label{3block-obj}
\min_{\m{u}} \; \frac{1}{2} \|\m{Fu} - \m{f}\|^2
+ \alpha \|\m{u}\|_{TV} + \beta \|\g{\Psi} \tr \m{u}\|_1,
\end{equation}
where $\m{f}$ is the given image data, $\m{F}$ is a matrix describing
the imaging device, $\|\cdot\|_{TV}$ is the total variation norm,
$\| \cdot \|_1$ is the $\ell_1$ norm,
$\g{\Psi}$ is a wavelet transform, and
$\alpha>0$ and $\beta >0$ are weights.
The first term in the objective is the data fidelity term,
while the next two terms are for regularization;
they are designed to enhance edges and increase image sparsity.
In our experiments, $\g{\Psi}$ is a normalized Haar wavelet 
with four levels and $\g{\Psi} \g{\Psi} \tr = I$.
The problem   (\ref{3block-obj}) is equivalent to 
\begin{equation}\label{3block-equiv}
\min_{(\m{u}, \m{w}, \m{z})} \; \frac{1}{2} \|\m{Fu} - \m{f}\|^2
+ \alpha \|\m{w}\|_{1,2} + \beta \|\m{z}\|_1 \; \mbox{subject to }
\m{Bu} = \m{w}, \; \g{\Psi} \tr \m{u} = \m{z},
\end{equation}
where $\m{Bu} = \nabla \m{u}$ and $(\nabla \m{u})_i$
is the vector of finite differences in the image along the
coordinate directions at the i-th pixel in the image, while
\[
\|\m{w}\|_{1,2} = \sum_{i=1}^N \| (\nabla \m{u})_i\|_2,
\]
where $N$ is the total number of pixels in the image.

The problem  (\ref{3block-equiv}) has the structure appearing in
(\ref{Prob})--(\ref{ProbM}) with
\[
\begin{array}{c}
\begin{array}{ll}
f_1 (\m{u})= 1/2  \|\m{Fu} - \m{f}\|^2, & h_1 := 0, \\[.05in]
f_2 := 0, & h_2(\m{w}) = \|\m{w}\|_{1,2}, \\[.05in]
f_3 := 0, & h_3(\m{z}) = \|\m{z}\|_1, \\[.05in]
\end{array} \\
\m{A}_1 = \left( 
\begin{array}{l}
 \m{B} \\
 \g{\Psi} \tr
\end{array}
 \right),
\quad
\m{A}_2 = \left( 
\begin{array}{r}
-\m{I} \\
 \m{0}
\end{array}
\right), \quad
\m{A}_3 =
\left( 
\begin{array}{r}
 \m{0} \\
-\m{I}
\end{array}
\right),
\quad \mbox{and} \quad
\m{b} =
\left( 
\begin{array}{r}
\m{0} \\
\m{0}
\end{array}
\right).
\end{array}
\]
When solving the test problems using accelerated BOSVS,
we use choose $\alpha^l$ and $\delta^l$ as in (\ref{AG_linesearch}).
Since $f_2 = f_3 = 0$, the line search condition holds automatically,
and the second and third subproblems are solved in closed form,
due to the simple structure of $h_2$ and $h_3$.
Only the first subproblem is solved inexactly.
At iteration $k$, the solution of this subproblem approximates the
solution of
\begin{eqnarray}\label{u-subprob}
 \min_{\m{u}} \; L_1^k (\m{u}) &: = &
 \frac{1}{2} \|\m{Au} - \m{b}\|^2 +
\frac{\rho}{2} \| \m{Bu} - \m{w}^k + \rho^{-1} \g{\lambda}^k\|^2 \\
&& \quad +  \frac{\rho}{2} \|  \g{\Psi} \tr \m{u} - \m{z}^k
+ \rho^{-1} \g{\mu}^k\|^2, \nonumber
\end{eqnarray}
where $\g{\lambda}^k$ and $\g{\mu}^k$ are the Lagrange multipliers
at iteration $k$ for the constraints
$\m{Bu} =$ $\m{w}$ and $\g{\Psi}\tr \m{u} =$ $\m{z}$ respectively.

The stopping condition for the inner loop of either multistep or
accelerated BOSVS required that $\Gamma_i^k \ge \Gamma_i^{k-1}$.
To improve efficiency, we replaced this condition by
\[
l_i^k \ge l_i^{k-1} \quad \mbox{or} \quad \Gamma_i^k \ge \Gamma_i^{k-1},
\]
where $l_i^k$ is the number of iterations performed
by the inner loop for block $i$ at iteration $k$.
For all the algorithm, we chose the
initial  $\delta_0^l$ in the line search using the BB approximation,
which is given in Step~1a of generalized BOSVS.
Moreover, when $\Gamma_i^k < \Gamma_i^{k-1}$, we increase $\delta_{\min,i}$
by setting
\[
\delta_{\min,i} := \tau \delta_{\min,i},
\]
where $\tau = 1.1$ in our numerical experiments.
When $\delta_{\min,i}$  is sufficiently large,
we have $\delta_0^l = \delta_{\min,i}$ and  the line search
condition in the algorithms is satisfied by $\delta_0^l$;
that is, $\delta^l= \delta_0^l =\delta_{\min,i}$.
Consequently, when $\delta_{\min,i}$  is sufficiently large, we have 
\[
\Gamma_i^k = \sum_{l=1}^{l_i^k} \frac{1}{\delta^l} =
\frac{ l_i^k}{\delta_{\min,i}},
\]
and the relaxed stopping condition $l_i^k \ge l_i^{k-1}$
implies that $\Gamma_i^k \ge \Gamma_i^{k-1}$, the original stopping condition.
Since $\tau > 1$, it follows that $\Gamma_i^k < \Gamma_i^{k-1}$ for only
a finite number of iterations, and hence, $\Gamma_i^k \ge \Gamma_i^{k-1}$ for
$k$ sufficiently large.
This ensures the global convergence of the algorithms.

Another improvement to efficiency was achieved by further relaxing the
line search criterion.
In particular, for the line search in generalized BOSVS (Step~1b),
we replaced the right side $f_i(\m{x}_i^{k+1})$ by
$f_i(\m{x}_i^{k+1}) - \epsilon^k$ where $\epsilon^k \ge 0$ is a summable
sequence.
In the line search of multistep BOSVS (Step~1b),
$f_i(\m{u}_i^l)$ was replaced by $f_i(\m{u}_i^l) - \pi^l$, where
$\pi^l = \epsilon^k \delta^l\omega^l$ with $\omega^l$ a summable sequence.
In the line search of accelerated BOSVS (Step~1a),
we replaced $f_i(\m{a}_i^l)$ by $f_i(\m{a}_i^l) - \pi^l$, where
$\pi^l = \epsilon^k \omega^l/\gamma^l$.
It can be proved that when the line search is relaxed
in this way using summable sequences, there is no
effect on the global convergence theory;
these $\epsilon^k$ and $\pi^l$
terms need to be inserted in each inequality in the analysis,
but in the end, the steps and the conclusions are unchanged.
On the other hand, when the line search is relaxed,
it can terminate sooner, and the algorithms can be more efficient.
For the numerical experiments, we took $\epsilon^k = 10/k^{1.1}$.
For multistep BOSVS, $\omega^l = 1/(\gamma^l)^{1.2}$,
while for accelerated BOSVS, $\omega^l = 1/(\gamma^l)^{0.6}$.
Since $\gamma^l = O(1/l)$ for multistep BOSVS and
$\gamma^l = O(1/l^2)$ for accelerated BOSVS, the $\omega^l$ sequences
are summable.

In all the algorithms, we use the following parameters:
\[ 
\delta_{\min} = 10^{-10},\; \delta_{\max} = 10^{10}, \; \alpha= 0.999,
\; \sigma=10^{-5},  \; \eta = 3, \; \mbox{and } \tau = 1.1.
\] 
For the inner loop stopping condition,
we took $\psi(t) = \min \{0.1 t, t^{1.1}\}$ in
multistep BOSVS, and $\psi(t) = 0.5t$ in accelerated BOSVS,
while in Step~2 of the ADMM template Algorithm~\ref{ADMMcommon},
we took $\theta_1 =$ $10^{-6} \sqrt{\rho}$,
$\theta_2 = \sqrt{\rho}$, and $\theta_3 = 10^{-6} \sqrt{\sigma/(1-\alpha)}$.
For comparison, we provide numerical results based on the algorithm
in \cite{HTXY2013} where we use MATLAB's
conjugate gradient routine {\sc cgs}
to solve the subproblem (\ref{u-subprob}) almost exactly,
stopping when $\|\nabla  L_1^k (\m{u})\| \le 10^{-6}$.
All the codes were implemented in MATLAB (version R2014a).
The following figures show the relative objective error
$(\Phi(\m{u}^k) - \Phi^*)/\Phi^*$  versus CPU time,
where $\Phi^*$ is the optimal function value of
(\ref{3block-obj}) obtained by applying accelerated BOSVS until
the eighth digit of the relative objective value did not change
in four consecutive iterations. 

\begin{figure}
\centering
 {\rm (a)} {\includegraphics[width=.45\textwidth]{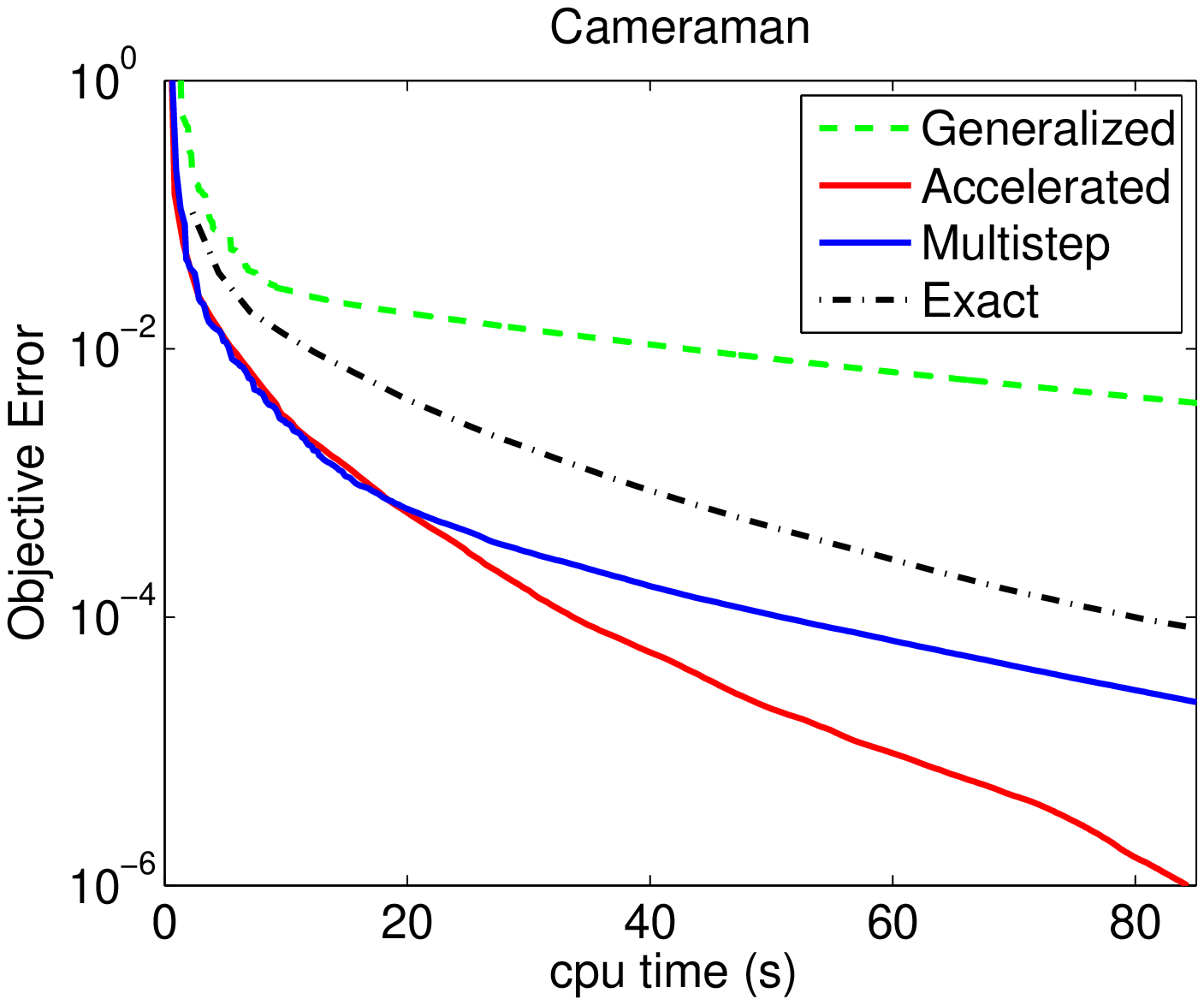}}
 {\rm (b)} {\includegraphics[width=.45\textwidth]{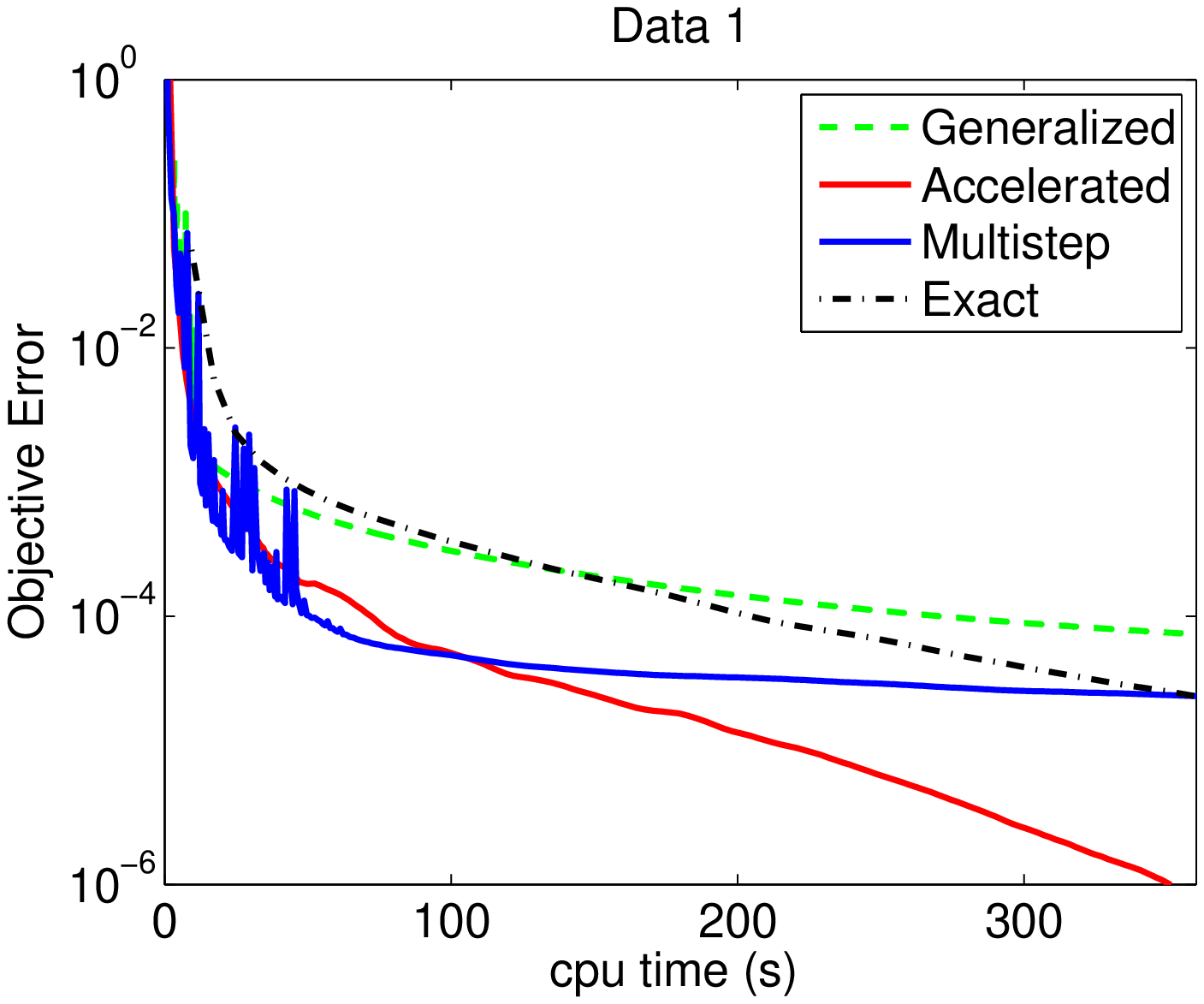}}
\\[.2in]
 {\rm (c)} {\includegraphics[width=.45\textwidth]{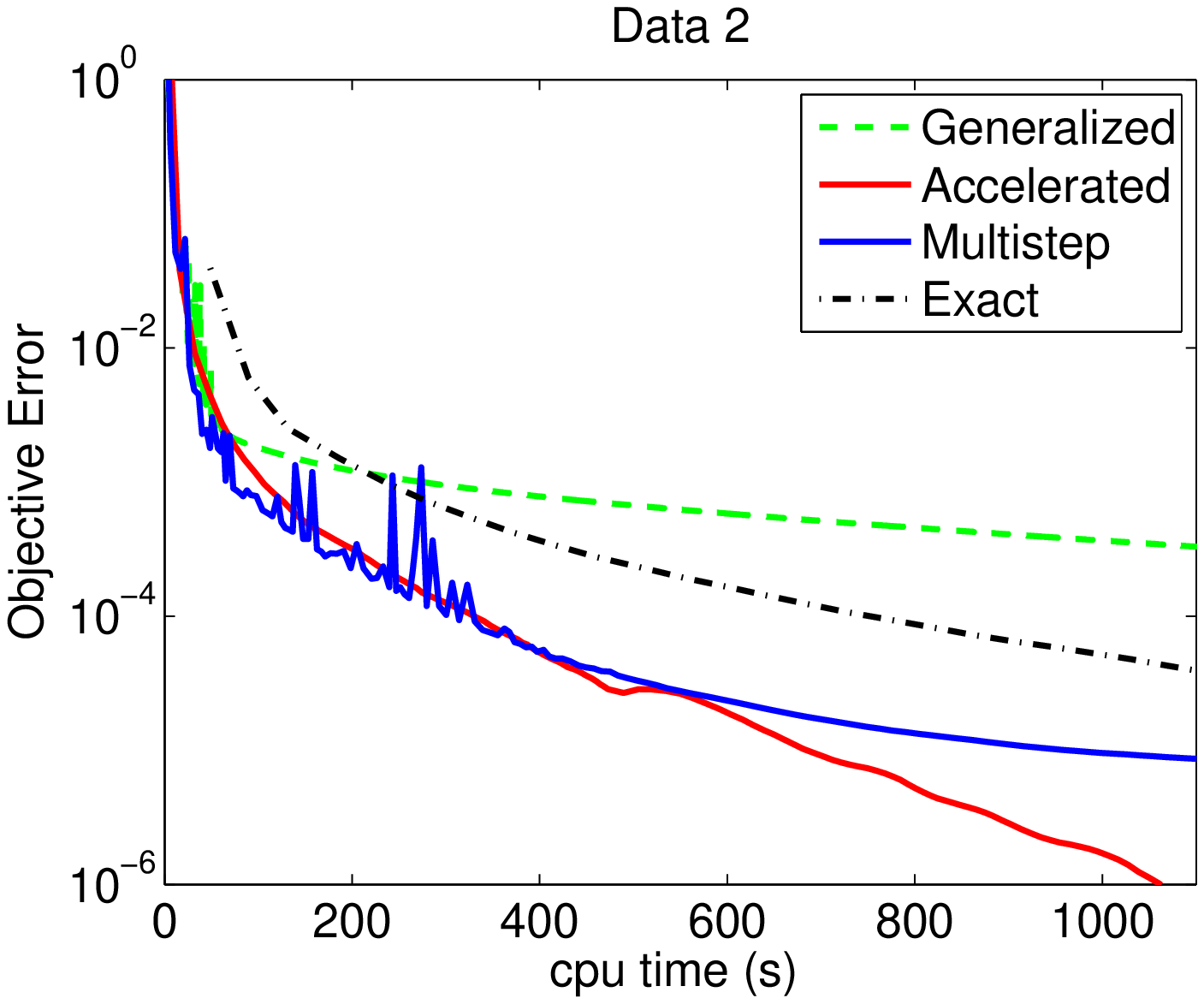}}
 {\rm (d)} {\includegraphics[width=.45\textwidth]{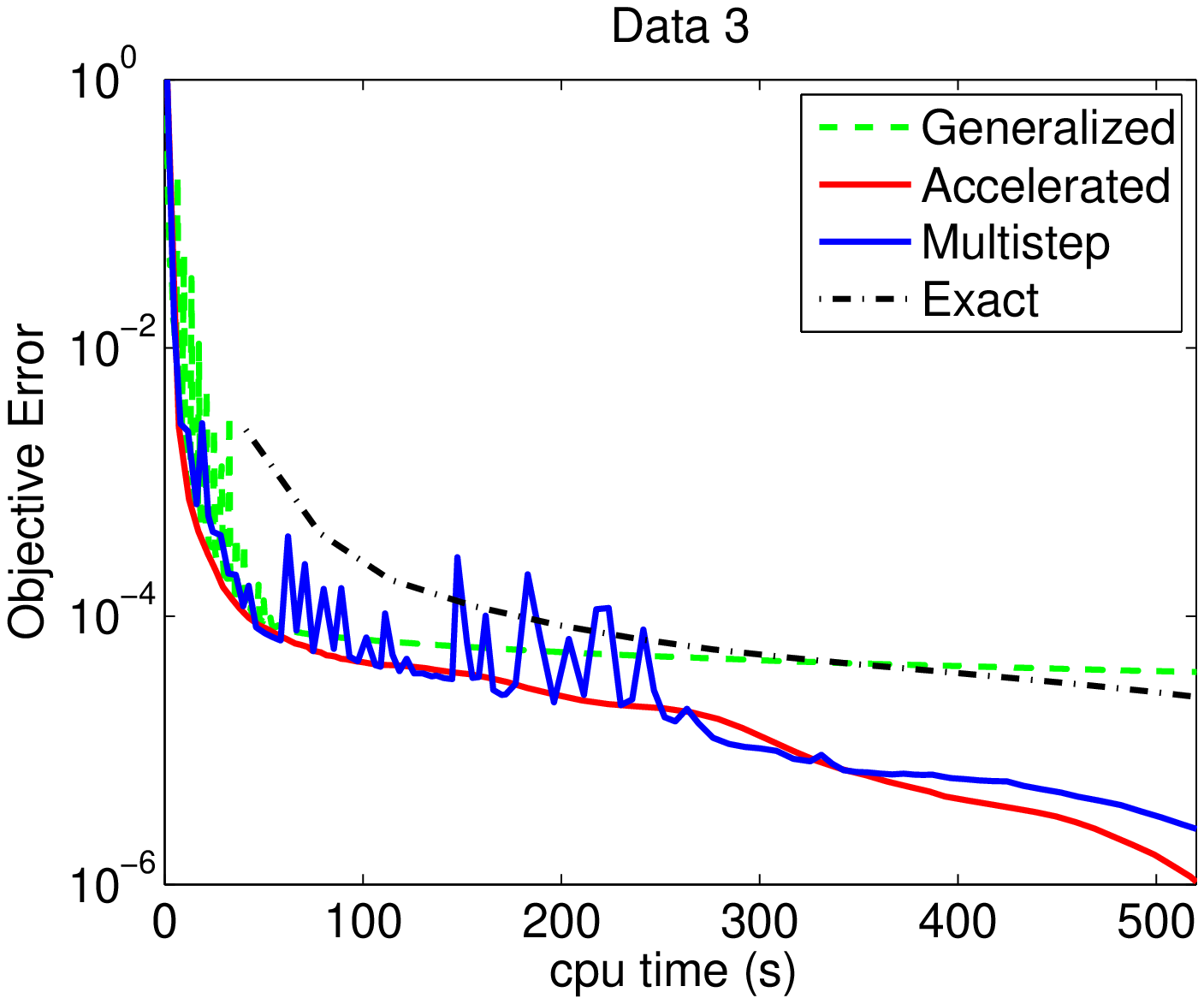}}
\caption{Base-10 logarithm of the relative objective error versus CPU time
for the test problems.
}
\label{error_plots}
\end{figure}

The first experiment employs an image deblurring problem from \cite{Afonso09b}.
The original image is the well-known Cameraman
image of size $256 \times 256$ and the observed data $\m{f}$ in
(\ref{3block-obj}) is a blurred image obtained by imposing a uniform blur of
size $9 \times 9$ with Gaussian noise and SNR of $40$dB.
The weights in (\ref{3block-obj}) are $\alpha = 0.005$ and $\beta = 0.001$,
and the penalty parameter $\rho = 5 \times 10^{-4}$.
Figure~\ref{error_plots}(a) shows the base-10 logarithm of the relative
objective error versus CPU time.
In this problem where the subproblems are
relatively easy, generalized BOSVS is significantly slower than
the exact, multistep, and accelerated algorithms,
while both multistep and accelerated BOSVS were faster
than the exact scheme.

The second set of test problems,
which arise in partially parallel imaging (PPI), are found in \cite{chy13}.
The observed data, corresponding to 3 different images, are
denoted data~1, data~2, and data~3.
For these test problems, the weights in (\ref{3block-obj}) are
$\alpha = 10^{-5}$ and $\beta = 10^{-6}$,
and the penalty parameter $\rho = 10^{-3}$.
The performance of the algorithms is shown in
Figure~\ref{error_plots}(b)--(d).
These test problems are much more difficult than the first
problem since $\m{F}$ is large, relatively dense, and ill conditioned.
In this case, all the inexact algorithms are faster
than the exact algorithm initially.
The exact algorithm becomes faster than generalized BOSVS when the
relative error is around $10^{-3}$ or $10^{-4}$.
Accelerated BOSVS is always significantly faster than the exact algorithm.

%
\section{Conclusion}
\label{conclusion}
Three inexact alternating direction multiplier methods were presented for
solving separable convex linearly constrained optimization problems,
where the objective function is the sum of smooth and
relatively simple nonsmooth terms.
The nonsmooth terms could be infinite, so the algorithms and analysis
included problems with additional convex constraints.
These algorithms all originate from the 2-block
variable stepsize BOSVS scheme of \cite{chy13, hyz14} which employs
indefinite proximal terms and linearized subproblems.
The 2-block scheme was generalized to a multiblock scheme using
a back substitution process to generate
an auxiliary sequence $\m{y}^k$ that played the role of
$\m{x}^k$ in the original, potentially divergent \cite{chyy2016},
multiblock ADMM (\ref{adm}).
The three new methods, called generalized, multistep, and accelerated BOSVS,
correspond to different accuracy levels when solving the ADMM subproblems.
Generalized BOSVS employed only one iteration in the subproblems,
while multistep and accelerated BOSVS performed multiple iterations until
the iteration change was sufficiently small.
The multistep and accelerated schemes differed in the rate with which
they solved the the subproblems.
If $l$ was the number of iterations in the subproblem, then multistep BOSVS
had a convergence rate of $O(1/l)$, while accelerated BOSVS had a
convergence rate of $O(1/l^2)$.
Global convergence was established for all the methods.
Numerical experiments were performed using image reconstruction problems.
The accelerated BOSVS algorithm had the best performance when compared
with either the other inexact algorithms, or the
exact algorithm of \cite{HTXY2013}.
\bibliographystyle{siam}

\end{document}